\documentclass[11pt]{amsart}
  \usepackage{amsmath,amssymb, graphics, amscd,latexsym}

  \makeatletter

  \newtheorem{Theorem}{Theorem}
  \newtheorem{Lemma}[Theorem]{Lemma}
  \newtheorem{Corollary}[Theorem]{Corollary}
  \newtheorem{Proposition}[Theorem]{Proposition}

  \theoremstyle{remark}
  \newtheorem{Remark}[Theorem]{Remark}
  
  \newtheorem{Example}[Theorem]{Example}

  \newtheorem{Assertion}[Theorem]{Assertion}

\begin{document}
\newcommand{\eps}{\varepsilon}
\newcommand{\om}{\omega}
\newcommand\Om{\Omega}
\newcommand\la{\lambda}
\newcommand\vphi{\varphi}
\newcommand\vrho{\varrho}
\newcommand\al{\alpha}
\newcommand\La{\Lambda}
\newcommand\si{\sigma}
\newcommand\be{\beta}
\newcommand\Si{\Sigma}
\newcommand\ga{\gamma}
\newcommand\Ga{\Gamma}
\newcommand\de{\delta}
\newcommand\De{\Delta}

\newcommand\cA{\mathcal  A}
\newcommand\cB{\mathcal B}
\newcommand\cD{\mathcal  D}
\newcommand\cM{\mathcal  M}
\newcommand\cN{\mathcal  N}
\newcommand\cT{\mathcal  T}
\newcommand\cP{\mathcal  P}
\newcommand\cp{\mathcal  p}
\newcommand\cQ{\mathcal  Q}
\newcommand\cG{\mathcal G}
\newcommand\cq{\mathcal  q}
\newcommand\cc{\mathcal  c}
\newcommand\cs{\mathcal  s}
\newcommand\cS{\mathcal  S}
\newcommand\ct{\mathcal  t}
\newcommand\cZ{\mathcal  Z}
\newcommand\cR{\mathcal  R}
\newcommand\cu{\mathcal  u}
\newcommand\cU{\mathcal  U}
\newcommand\cI{\mathcal  I}
\newcommand\cJ{\mathcal  J}
\newcommand\co{\mathcal  o}
\newcommand\cO{\mathcal  O}
\newcommand\cv{\mathcal  v}
\newcommand\cV{\mathcal  V}
\newcommand\cx{\mathcal  x}
\newcommand\cX{\mathcal  X}
\newcommand\cw{\mathcal  w}
\newcommand\ck{\mathcal  k}
\newcommand\cK{\mathcal  K}
\newcommand\cW{\mathcal  W}
\newcommand\cz{\mathcal  z}
\newcommand\cy{\mathcal  y}
\newcommand\ca{\mathcal  a}
\newcommand\ch{\mathcal  h}
\newcommand\cH{\mathcal  H}
\newcommand\cF{\mathcal F}
\newcommand\bfG{\mbox {\bf  G}}
\newcommand\bfg{\mbox {\bf  g}}
\newcommand\bfC{\mbox {\bf  C}}
\newcommand\bfN{\mbox {\bf  N}}
\newcommand\bfT{\mbox {\bf  T}}
\newcommand\bfP{\mbox {\bf  P}}
\newcommand\bfp{\mbox {\bf  p}}
\newcommand\bfQ{\mbox {\bf  Q}}
\newcommand\bfq{\mbox {\bf  q}}
\newcommand\bfc{\mbox {\bf  c}}
\newcommand\bfs{\mbox {\bf  s}}
\newcommand\bfS{\mbox {\bf  S}}
\newcommand\bft{\mbox {\bf  t}}
\newcommand\bfZ{\mbox {\bf  Z}}
\newcommand\bfR{\mbox {\bf  R}}
\newcommand\bfu{\mbox {\bf  u}}
\newcommand\bfU{\mbox {\bf  U}}
\newcommand\bfo{\mbox {\bf  o}}
\newcommand\bfO{\mbox {\bf  O}}
\newcommand\bfv{\mbox {\bf  v}}
\newcommand\bfV{\mbox {\bf  V}}
\newcommand\bfx{\mbox {\bf  x}}
\newcommand\bfX{\mbox {\bf  X}}
\newcommand\bfw{\mbox {\bf  w}}
\newcommand\bfk{\mbox {\bf  k}}
\newcommand\bfK{\mbox {\bf  K}}
\newcommand\bfW{\mbox {\bf  W}}
\newcommand\bfz{\mbox {\bf  z}}
\newcommand\bfy{\mbox {\bf  y}}
\newcommand\bfa{\mbox {\bf  a}}
\newcommand\bfh{\mbox {\bf  h}}
\newcommand\bfH{\mbox {\bf  H}}
\newcommand\bfJ{\mbox {\bf  J}}
\newcommand\bfj{\mbox {\bf  j}}
\newcommand\bbC{\mbox {\mathbb C}}
\newcommand\bbN{\mbox {\mathbb N}}
\newcommand\bbT{\mbox {\mathbb T}}
\newcommand\bbP{\mbox {\mathbb P}}
\newcommand\bbQ{\mbox {\mathbb Q}}
\newcommand\bbS{\mbox {\mathbb S}}
\newcommand\bbZ{\mbox {\mathbb Z}}
\newcommand\bbR{\mbox {\mathbb R}}
\newcommand\bbU{\mbox {\mathbb U}}
\newcommand\bbO{\mbox {\mathbb O}}
\newcommand\bbV{\mbox {\mathbb V}}
\newcommand\bbX{\mbox {\mathbb X}}
\newcommand\bbK{\mbox {\mathbb K}}
\newcommand\bbW{\mbox {\mathbb W}}
\newcommand\bbH{\mbox {\mathbb H}}

\newcommand\apeq{\fallingdotseq}
\newcommand\Lrarrow{\Leftrightarrow}
\newcommand\bij{\leftrightarrow}
\newcommand\Rarrow{\Rightarrow}
\newcommand\Larrow{\Leftarrow}
\newcommand\nin{\noindent}
\newcommand\ninpar{\par \noindent}
\newcommand\nlind{\nl \indent}
\newcommand\nl{\newline}
\newcommand\what{\widehat}
\newcommand\tl{\tilde}
\newcommand\wtl{\widetilde}
\newcommand\order{\mbox{\text{order}\/}}
\newcommand\GL{\text{GL}\/}
\newcommand\PGL{\text{PGL}\/}
\newcommand\Spec{\text{Spec}\/}
\newcommand\weight{\text{weight}\/}
\newcommand\ord{\text{ord}\/}
\newcommand\Int{\text{Int}\/}
\newcommand\grad{\text{grad}\/}
\newcommand\Ind{\text{Ind}\/}
\newcommand\Disc{\text{Disc}\/}
\newcommand\Ker{\text{Ker}\/}
\newcommand\Image{\text{Image}\/}
\newcommand\Coker{\text{Coker}\/}
\newcommand\Id{\text{Id}\/}
\newcommand\id{\text{id}}
\newcommand\dsum{\text{\amalg}}
\newcommand\val{\text{val}}
\newcommand\mv{\text{m-vector}}
\newcommand\iv{\text{i-vector}}
\newcommand\minimum{\text{minimum}\/}
\newcommand\modulo{\text{modulo}\/}
\newcommand\Aut{\text{Aut}\/}
\newcommand\PSL{\text{PSL}\/}
\newcommand\Res{\text{Res}\/}
\newcommand\rank{\text{rank}\/}
\newcommand\codim{\text{codim}\/}
\newcommand\msdim{\text{ms-dim}\/}
\newcommand\icodim{\text{i-codim}\/}
\newcommand\ocodim{\text{o-codim}\/}
\newcommand\emdim{\text{ems-dim}\/}

\newcommand\Cone{\text{Cone}\/}
\newcommand\maximum{\text{maximum}\/}
\newcommand\Vol{\text{Vol}\/}
\newcommand\Coeff{\text{Coeff}\/}
\newcommand\lcm{\text{lcm}\/}
\newcommand\degree{\text{degree}\/}
\newcommand\Pol{\cal {POL}}
\newcommand\Ts{Tschirnhausen}
\newcommand\TS{Tschirnhausen approximate}
\newcommand\Stab{\text{Stab}\/}
\newcommand\civ{complete intersection variety}
\newcommand\nipar{\par \noindent}
\newcommand\wsim{\overset{w}{\sim}}
\newcommand\Gr{\bfZ_2*\bfZ_3}
\newcommand\QED{~~Q.E.D.}
\newcommand\bsq{$\blacksquare$}
\newcommand\bff{\mbox {\bf  f}}
\newcommand\newcommandby{:=}
\newcommand\inv{^{-1}}
\newcommand\nnt{(\text{nn-terms})}
\renewcommand{\subjclassname}{\textup{2000} Mathematics Subject Classification}
\def\mapright#1{\smash{\mathop{\longrightarrow}\limits^{{#1}}}}
\def\mapleft#1{\smash{\mathop{\longleftarrow}\limits^{{#1}}}}
\def\mapdown#1{\Big\downarrow\rlap{$\vcenter{\hbox{$#1$}}$}}
\def\mapdownn#1#2{\llap{$\vcenter{\hbox{$#1$}}$}\Big\downarrow\rlap{$\vcenter{\hbox{$#2$}}$}}
\def\mapup#1{\Big\uparrow\rlap{$\vcenter{\hbox{$#1$}}$}}
\def\mapupp#1#2{\llap{$\vcenter{\hbox{$#1$}}$}\Big\uparrow\rlap{$\vcenter{\hbox{$#2$}}$}}

\def\rdown#1{\searrow\rlap{$\vcenter{\hbox{$#1$}}$}}
\def\semap#1{\searrow\rlap{$\vcenter{\hbox{$#1$}}$}}

\def\rup#1{\nearrow\rlap{$\vcenter{\hbox{$#1$}}$}}
\def\nemap#1{\nearrow\rlap{$\vcenter{\hbox{$#1$}}$}}
\def\ldown#1{\swarrow\rlap{$\vcenter{\hbox{$#1$}}$}}
\def\swmap#1{\swarrow\rlap{$\vcenter{\hbox{$#1$}}$}}
\def\lup#1{\nwarrow\rlap{$\vcenter{\hbox{$#1$}}$}}
\def\nwmap#1{\nwarrow\rlap{$\vcenter{\hbox{$#1$}}$}}
\def\defby{:=}
\def\eqby#1{\overset {#1}\to =}
\def\inv{^{-1}}
\def\bnu{{(\nu)}}
\def\ocup#1{\underset{#1}\cup}
  \title[Geometry of reduced sextics of  torus type
  ]
  {Geometry of reduced sextics of  torus type}

  \author
  [M. Oka ]
  {Mutsuo Oka }
  \address{{Mutsuo Oka}\\
  {Department of Mathematics}\\
  {Tokyo Metropolitan University}\\
  {1-1 Mimami-Ohsawa, Hachioji-shi}\\
  {Tokyo 192-0397}\\
  {\rm{E-mail}}: {\rm oka@comp.metro-u.ac.jp}}


  \begin{abstract}
  In \cite{Oka-Pho2}, we gave a classification of the 
 configurations  of
   irreducible sextic of torus type.
  In this paper, we give a
  classification of the singularities on reducible sextics of  torus
   type. We determine the components types and  the geometry of the
   components
for each configurations.
  \end{abstract}
  \maketitle

  \pagestyle{headings}
  \section{Introduction and statement of the result}
  In our previous paper \cite{Oka-Pho2}, we have classified the
  configurations of the singularities on irreducible sextics of torus type.
  In this note, we classify the configuration of the singularities
   of reducible sextics of
  torus type. 
  We use the same notations as in \cite{Oka-Pho2}.

   We denote by $\Si_{in}$  and $\Si_{out}$ the
   configurations
  of the inner singularities and of the outer
   singularities respectively.
For the classification of the configurations
of the reduced sextics of torus type,
it is less important to 
distinguish inner or outer singularities  but what is more important
is to know  the singularities of  the irreducible components
and their intersections.
We put $\Si_{red}=\Si_{in}\cup\Si_{out}$, and we call it {\em the reduced
   configuration}.
Let $B_{i_1},\dots, B_{i_k}$ be the irreducible components of $C$.
We call $\{\deg B_{i_1},\dots \deg B_{i_k}\}$
{\em  the  component type}  of a reducible sextics $C$.
In this note, we assume that  the curves $B_i,B_i',\dots$   are 
   irreducible
and their degrees are the same with the indices. Thus, for example,
  $C=B_1+B_1'+B_4$
   implies 
 that 
$C$ has three  components of degree 1, 1, 4.
  The configurations of the singularities of $B_{i}$ is  denoted by $\Si(B_{i})$.
  We say that  $C$ 
has the maximal rank if
$C$ has only simple singularities and the total Milnor number is 19.
  We denote  configurations with maximal rank by upper suffix $mr$, like
 $[A_{11},2A_2, D_4]^{mr}$.
The classification of reduced sextics of torus type with 
only simple singularities is given in Theorem 1 and the classification
 for
the  other case is given in Theorem 2.

\subsection{Reduced sextics with simple singularities}
 We first classify the reduced sextics with simple singularities.
  \begin{Theorem}\label{Reduced-Simple}
  The classifications of singularities on reducible sextics
  of torus type with only simple singularities are given  as follows.

  \begin{enumerate}
  \item  $\Si_{in}=[A_5,4A_2]:$  $C=B_5+B_1$ and 
 $[A_5,  4A_2, 2A_1]_2,  [A_5, 4A_2, 3A_1]_2$,
$[A_5, 5A_2, 2A_1]_2$, $[A_5, 4A_2, A_3, 2A_1]_2$,
$[A_5, 4A_2, 4A_1]$, 
$[A_5, 4A_2, A_3]_2$,$[A_5, 4A_2,A_3, A_1]_2,[A_5, 4A_2,A_3, 2A_1]_3$,
$[A_5,4A_2,D_4]_2, [A_5,4A_2,D_5]_2$.

  \item $\Si_{in}=[2A_5,2A_2]${\rm :}
\begin{enumerate}
\item $C=B_1+B_5${\rm :}
 $[2A_5,2A_2,2A_1]_2,  [2A_5,2A_2,3A_1]_1,  [2A_5,3A_2,2A_1]_1,
[2A_5,2A_2,A_3]_2$, $[2A_5,2A_2,A_3,A_1]_1$,
$[2A_5,2A_2,D_4]_1$, $[2A_5,2A_2,D_5]_1^{mr}$.
\item $C= B_1+B_1'+\, B_4${\rm :}
$[2A_5,2A_2,3A_1]_2,  [2A_5,2A_2,4A_1]$, $[2A_5,2A_2,D_4]_2$.
\item $C= B_2+B_4${\rm :}
 $ [2A_5,2A_2,2A_1]_3$, $[2A_5,2A_2,3A_1]_3$, $[2A_5, 3A_2, 2A_1]_2$,
 $[2A_5,2A_2,A_3]_3$, $[2A_5,2A_2,A_3,A_1]_2$,
$[2A_5,2A_2,D_4]_3,\, [2A_5,2A_2,D_5]_2^{mr}$.
\item $C=B_3+B_3'${\rm :}
$ [2A_5,2A_2,3A_1]_4,\, [2A_5,2A_2,A_3,A_1]_3,
  [ 3A_5,2A_2]^{mr}$.
\end{enumerate}
  \item $\Si_{in}=[E_6,A_5,2A_2]${\rm :} $C=B_1+B_5$ and
$ [E_6, A_5, 2A_2, 2A_1]_2$, $[E_6, A_5, 2A_2, 3A_1]$,
$[E_6, A_5, 2A_2, A_3]_2$, $[E_6, A_5, 2A_2, A_3,A_1]^{mr}$.
  \item  $\Si_{in}=[3A_5]${\rm :}
\begin{enumerate}
\item $C=B_1+B_5${\rm :}
  $[3A_5,  2A_1]_1$, $[3 A_5,A_3]_1$.
\item $C=B_2+B_4${\rm :} $[3A_5,2A_1]_2$, $[3A_5,A_3]_2$.
\item $C=B_1+B_1'+B_4${\rm :}
  $[3A_5, 3A_1]_1$, $[3A_5, D_4]_1^{mr}$.
\item $C=B_3+B_3'${\rm :}  $[3A_5]_2$, $[3A_5, A_1]_2$,
$[3A_5,A_2]_2$, $[3A_5, 2A_1]_3$, $[3A_5, A_1, A_2]$, $[3A_5, 2A_2]^{mr}$.
\item $C=B_1+B_2+B_3${\rm :}
$ [3 A_5, 2A_1],  [3 A_5,  3A_1]_1, 
   [3A_5,A_2,2A_1]^{mr},  [3A_5,\,A_3]_3$, $[3A_5,\,A_3, A_1]^{mr}$.
\item $C=B_1+B_1'+B_1''+B_3${\rm :}  $[3A_5,3A_1]_2$, $[3A_5,4A_1]^{mr}$,
$[3A_5,D_4]_2^{mr}$.
\item $C=B_2+B_2'+B_2''${\rm :} $[3A_5,3A_1]_3$, $[3A_5,D_4]_3^{mr}$.
\end{enumerate}
  \item $\Si_{in}=[2A_5,E_6]${\rm :}
\begin{enumerate}
\item $C=B_1+B_5${\rm :}
$  [E_6,2A_5,2A_1]_1$, $[E_6,2A_5,A_3]_1^{mr}$.
\item $C=B_2+B_4${\rm :}
  $[E_6,2A_5,2A_1]_2$, $[E_6,2A_5,A_3]_2^{mr}$.
\item $C=B_1+B_1'+B_4${\rm :}
$[E_6,2A_5,3A_1]^{mr}$.
\end{enumerate}
  \item $\Si_{in}=[A_8,A_5,A_2]${\rm :}  $C=B_1+B_5$ and 
$  [A_8, A_5,A_2, 2A_1]_2$, $[A_8, A_5,A_2,  3A_1]$,
$[A_8,A_5,A_2, A_3]_2$,   $[A_8,A_5,A_2,A_3,A_1]^{mr}$,
$ [A_8,A_5,A_2,D_4]^{mr}$.
  \item $\Si_{in}=[A_{11},2A_2]${\rm :}
\begin{enumerate}
\item $C=B_2+B_4${\rm :}
 $[A_{11},2A_2,2A_1]_2$, $[A_{11},2A_2, 3A_1]_1$, 
  $[A_{11},3A_2,2A_1]^{mr}$, $[A_{11},2A_2,A_3]_2$,
$[A_{11},2A_2,D_4]^{mr}$.
\item $C=B_3+B_3'${\rm :}
$[A_{11},2A_2,3A_1]_2$, $[A_{11},2A_2,A_3,A_1]$.
\end{enumerate}
  \item $\Si_{in}=[A_{11},A_5]${\rm :}
\begin{enumerate}
\item  $C=B_1+B_5${\rm :}
$[A_{11},A_5, 2A_1]_1$, $[A_{11},A_5, A_3]_1^{mr}$.
\item $C=B_2+B_4${\rm :}
$[A_{11},A_5,2A_1]_2 ,  [A_{11},A_5, A_3]_2^{mr}$.
\item $C=B_3+B_3'${\rm :}
$[A_{11},A_5]_2$, $[A_{11},A_5,A_1]_2$, $[A_{11},A_5,A_2]_2$,
$[A_{11},A_5,2A_1]_3$, $[A_{11},A_5,A_2,A_1]^{mr}$.

\item $C=B_1+B_2+B_3${\rm :}
$[A_{11},A_5,2A_1]$, $[A_{11},A_5,3A_1]^{mr}$,
$[A_{11},A_5,A_3]_3^{mr}$.
\end{enumerate}
  \item $\Si_{in}=[A_{17}]${\rm :}  $C=B_3+B_3'$,
$[A_{17}]_2$, $[A_{17},A_1]_2$,  
$[A_{17},2A_1]^{mr}$,  $[A_{17},A_2]_2^{mr}$.
\end{enumerate}
  \end{Theorem}
Further geometrical informations are  explained in the proof in Section 3.
In the above theorem, the  lower index, like $[A_5, 4A_2, A_3, 2A_1]_2$,
is to distinguish other component with the same 
weak Zariski configuration. The index 1 is reserved for the irreducible
case if there exists an irreducible sextics. Thus the configuration
which start from index 2,  like $[A_5, 4A_2, A_3, 2A_1]_2$, implies that
	there is an irreducible sextics with the same configuration.
  \begin{Theorem}\label{Reduced-NonSimple}
  The reduced configurations with at least one non-simple singularities
  are given by the following.
  \begin{enumerate}
  \item  $B_{3,6}\in \Si_{in}:$
\begin{enumerate}
\item $C=B_1+B_5:$
$[B_{3,6},3A_2,A_1]_2$,  $[B_{3,6},3A_2,2A_1]$, 
 $[ B_{3,6},4A_2,A_1] $, $[B_{3,6},A_5,A_2,A_1]_1$,
$[B_{3,6},E_6,A_2,A_1] $, 
  $[B_{3,6},A_8,A_1] $.
\item $C=B_1+B_1'+B_4:$
 $[B_{3,6},A_5,A_2,2A_1]_1 $.
\item $C=B_1+B_2+B_3:$ $[B_{3,6},A_5,A_2,2A_1]_2 $.
 \item $C=B_2+B_4:$
$[B_{3,6},A_5,A_2,A_1]_2$.
\item $C=B_2+B_2'+B_2'':$
	$[2B_{3,6}]$.
\end{enumerate}
  \item $C_{3,7}\in \Si_{in}:$
\begin{enumerate}
\item $C=B_1+B_5:$
  $[C_{3,7},3A_2,A_1]_2$,  $[C_{3,7},3A_2,2A_1]$,  
$[C_{3,7},A_5,A_2,A_1]_1$,   $[C_{3,7},E_6,A_2,A_1]$,   $[C_{3,7},A_8,A_1]$.
\item $C=B_2+B_4:$ $[A_2,A_5,C_{3,7},A_1]_2$.
\item $C=B_1+B_1'+B_4:$ $[C_{3,7},A_5,A_2,2A_1]$.
\end{enumerate}
  \item $C_{3,8}\in \Si_{in}:$
\begin{enumerate}
\item {\rm (a-1)} $C=B_1+B_5$ and $(B_5,O)=A_3$:
$[C_{3,8},3A_2]_2$,  $[C_{3,8},3A_2,A_1]_1 $, $[C_{3,8},A_5,A_2]$,  
$[C_{3,8},E_6,A_2] $, $[C_{3,8},A_8] $. 

{\rm (a-2)}  $C=B_1+B_5$ and $(B_5,O)=A_5 :$
$[C_{3,8},3A_2,A_1]_2$.
\item $C=B_2+B_4:$ $[C_{3,8},A_5,A_2]$.
\item
 $C=B_1+B_2+B_3:$ $[C_{3,8},A_5,A_2,A_1]_1$,
$[C_{3,8},A_5,A_2, A_1]_2$.
\end{enumerate}
\item $C_{3,9}\in \Si_{in}:$
  $C=B_1+B_5$ and $[C_{3,9},2A_2,A_1]_2$, $[C_{3,9},2A_2,2A_1]$,
$[C_{3,9},3A_2,A_1]$, $[C_{3,9},A_5,A_1]$,  $[C_{3,9},E_6,A_1]$.
  \item $C_{3,12}\in \Si_{in}:$
\begin{enumerate}
\item $C=B_1+B_5:$ $[C_{3,12},A_2,A_1]_1$.
\item $C=B_2+B_4:$ $[C_{3,12},A_2,A_1]_2$.
\item $C=B_1+B_2+B_3:$ $[C_{3,12},A_2,2A_1]$.
\end{enumerate}
  \item $C_{6,6}\in \Si_{in}:$
\begin{enumerate}
\item $C=B_1+B_5:$ $[C_{6,6},2A_2,A_1]$.
\item  $C=B_1+B_1'+B_4:$ $[C_{6,6},2A_2,2A_1]_1$.
\item $C=B_3+B_3'$: $[C_{6,6},A_5]$.
\item $C=B_3+B_2+B_1$: $[C_{6,6},A_5,A_1]$
\item  $C=B_1+B_1'+B_2+B_2':$ $[C_{6,6},A_5,2A_1]_2$.
\end{enumerate}
  \item $C_{6,9}\in \Si_{in}:$ $C=B_1+B_5$,
$[C_{6,9},A_2,A_1]$.
  \item $C_{3,15}\in \Si_{in}:$ $C=B_1+B_5$, $[C_{3,15},A_1]$.
\item $B_{3,12}\in \Si_{in}$: $C=B_2+B_2'+B_2''$, $[B_{3,12}]$.
  \item $C_{6,12}\in\Si_{in} :$
\begin{enumerate}
\item $C=B_3+B_3':$
  $[C_{6,12}]$.
\item $C=B_1+B_2+B_3:$
$[C_{6,12},A_1]$.
\end{enumerate}
  \item $B_{4,6}\in \Si_{in}:$ $C=B_3+B_3'$ and
$[B_{4,6},A_5]$.
  \item $D_{47}\in \Si_{in}:$
\begin{enumerate}
\item  $C=B_1+B_5:$
$[D_{4,7},2A_2]$.
\item $C=B_1+B_2+B_3:$ $[D_{4,7},A_5]$.
\end{enumerate}
  \item $Sp_2\in \Si_{in}:$ $C=B_3+B_3'$ and $[Sp_2]$.
  \item $B_{6,6}\in \Si_{in}:$ $[B_{6,6}]$.
  \end{enumerate}
  \end{Theorem}
  \section{Preliminaries}
\subsection{Genus formula and the class formula}
Let  $C$ be an irreducible curve of a
  given degree $d$. Then the genus formula  is given as
  \begin{multline*}
   g(C)=\frac{(d-1)(d-2)}2 -\sum_{P\in \Si(C)}\de(P)\ge 0,\quad
  \de(P)=\frac{\mu(C,P)+r(C,P)-1}2
  \end{multline*}
  where $\mu(C,P)$ and $r(C,P)$is the Milnor number and the number of
  local irreducible components (\cite{Milnor}).
  Let $\de^*(C)=\sum_{P\in \Si(C)}\de(P)$.
  Using this criterion, we have
$ \de^*(C)\le 6, 3, 1, 0$  respectively for 
$d=5,\, 4,\, 3,\, 2$ for an irreducible curve $C$.
Now assume that $C$ is not irreducible.
  Let $C=B_{i_1}+\cdots B_{i_k}$ be the irreducible decomposition of $C$
  with $\degree(B_{i_k})=i_k$. We define $\de^*(C)=\sum_{j=1}^k\de^*(B_{i_j})$.
  \begin{Corollary}\label{sum-delta}
  Assume that 
  $C=B_{i_1}+\cdots B_{i_k}$ is the irreducible decomposition and $k\ge 2$.
  Then 
  \[
	\de^*(C)\le
  \begin{cases}
  6\quad &C=B_5+B_1\\
  3\quad & C=B_4+B_2,\, B_4+B_1+B_1'\\
  2\quad &C=B_3+B_3'\, \\
  1\quad & C=B_3+B_2+B_1,B_3+B_1+B_1'+B_1''\\
  0\quad & otherwise
  \end{cases}
  \]
  \end{Corollary}
The class formula describes the degree $n^*(C)$ of the dual curve
and it is given by the following formula(\cite{NambaBook}).
\[
 n^*(C)=d(d-1)-\sum_{P\in \Si(C)} (\mu(C,P)+m(C,P)-1)
\]
The number of flex points $i(C)$ counted with multiplicity is given by
\[
 i(C)=3d(d-2)-\sum_{P\in\Si(C)} \text{flex defext(C,P)}
\]
For the definition of flex defect, we refer Oka \cite{Oka-sextics}
  \subsection{Intersection singularities}
  Let $C$ be a plane curve and let $C^1,\dots, C^k$ be the irreducible 
components.
Let $P$ be  a singular point of $C$.
We say that $P$ is {\em a proper singularity} if 
$P\in C^i-\cup_{j\ne i}C^j$ for some component $C^i$. Otherwise
we say that  $(C,P)$  {\em an intersection singularity}.
  Assume that $C=C^1\cup C^2$, for example,
and  $P\in C^1\cap C^2$ and   $C^1,\, C^2 $ are non-singular at $P$ and let $\iota$ be
  the local intersection number. Then $(C,P)\cong A_{2\iota-1}$.
  Assume further that $C^1$ is a line and $\iota\ge 3$. Then we say that
  $C^1$
  is a flex tangent line of $C^2$.
  \begin{Proposition}
  Assume that $C$ is defined by $f(x,y)=0$ and 
  assume that the Newton boundary is non-degenerate.
  Let $\De_1,\dots, \De_k$ be the faces of $\Ga(f)$
  and assume that 
  $f_{\De_i}(x,y)=\prod_{j=1}^{\nu_i}(y^{a_i}-\al_j x^{b_i})$
  with $\gcd(a_i,b_i)=1$ and
  $\al_1,\dots,\al_{\nu_i}$ are mutually distinct.
  Then $C$ has  $\sum_{i=1}^k \nu_i$  local irreducible components
  which are defined 
  $(y^{a_i}-\al_j x^{b_i})+\text{(higher terms)}=0$. \end{Proposition}
  See for example \cite{Okabook}.

\begin{Example}\label{example-inter}
 1. 
Consider $D_4: y^2\,x +x^3=0$. Then $D_4$  can be an intersection
   singularity of three smooth components, $x=0, y\pm x=0$
  where each two of them intersect transversely.
  Similarly $D_5: y^2 x+x^4=0$ can be interpreted as an intersection 
  singularity of a line $x=0$ and a cusp $y^2+x^3=0$.
\nl
 2. Consider the singularity
  $C_{3,p}: ~y^3+y^2\,x^2-x^p=0$.
\nl
Case 1. Assume that $p$ is odd. Then $C_{3,p}$ has two local irreducible
components. One component is  smooth and defined by 
  $L: y+x^2+\text{(higher terms)}=0$ and 
  another  component $M$ is 
  defined
  by $ y^2-x^{p-2}+\text{(higher terms)}=0$  and it is an
  $A_{p-3}$-singularity and $I(L,M; O)=4$.
\nl
Case 2. Assume that $p$ is even and put $p=2m, \, m\ge 4$. Then $C_{3,2m}$ has
three
smooth components $L_1,L_2,L_3 $ where 
$L_1: y+x^2+\text{(higher terms)}=0$ and 
$L_2,L_3$ are defined by $ y\pm x^{m-1}+\text{(higher terms)}=0$. Note that 
$(L_2\cup L_3,O)\cong A_{2m-3}$ and $(L_1\cup L_2,O)\cong A_3$.
\end{Example}

We use the following notations
for the non-simple singularities  as in \cite{Pho}.
\begin{eqnarray*}
\begin{cases}
B_{p,q}:  \: y^p+x^q=0  \: (\text{Brieskorn-Pham type})\\
C_{p,q}:  \: y^p+x^q+x^2y^2=0,\quad
\frac 2p+\frac 2q\,<\, 1\\
D_{4,7}:  \: y^4+x^3y^2+a x^5 y+ bx^7=0,\quad a^2-4b\ne 0\\
Sp_1:  \: (y^2-x^3)^2+(xy)^3=0\\
Sp_2:  \: (y^2-x^3)^2-y^6=0
\end{cases}
\end{eqnarray*}
Hereafter  we only consider sextics of torus type
\[
 C:\quad f_2(x,\, y)^3+f_3(x,\, y)^2=0
\]
The notation $C_2: f_2(x,y)=0$ and $C_3: f_3(x,y)=0$ is used throughout
the paper.
  \subsection{Weak Zariski k-ple}
  A $k$-ple of reduced plane curves $\{C^1,\dots,C^k\}$ is called 
 {\em a weak Zariski $k$-ple}
  if $\degree(C^1)=\dots=\degree(C^k)$ and they have same reduced configuration of 
  singularities and 
  the topology of the complement $\bfP^2-C^j$ are all different.
We call $\Si(C^j)$ {\em a weak Zariski configuration}.
Note that $C^1,\dots, C^k$ may have different component types.
 Artal  has first   observed such a pair for $[A_{17}]$ or some others in
 sextics
\cite{Artal}. It is obvious that a Zariski pair is a weak Zariski
pair.

\subsection{Sextics of linear torus}
A sextics $C$  of torus type  is called {\em of linear torus type } if $C$ can be defined by
$ f(x,y)=f_2(x,y)^3+f_3(x,y)^2$  where $f_2(x,y)=(a x+by+c)^2$.
We may assume that $f_2=-y^2$ by a linear change of coordinates
so that $f$ is a product of linear forms
$f(x,y)=(f_3(x,y)+y^3)(f_3(x,y)-y^3)$.
It is easy to observe that the inner singularities
are on $y=f_3(x,0)=0$.
\begin{Proposition}
The possibility of inner configuration of sextics of linear torus type
is either $[3A_5]$, or $[A_{11},A_5]$ or $[A_{17}]$ for simple
 singularities
and for non-simple singularities, we have
$[C_{6,6},A_5]$, $[B_{4,6}, A_5]$, $[D_{4,7}, A_5]$,
$[C_{6,12}]$, $[Sp_2]$ and $[B_{6,6}]$.
\end{Proposition}
{\em Proof.}
Assume that $f_3(\alpha,0)=0$.  Assume that $\alpha$ is a
simple  solution of $f_3(x,0)=0$ (respectively
a solution  of multiplicity 2 or 3). Put $P=(\alpha,0).$ Then
$I(y^2,C_3;P)=2$ (resp. 4 or 6).
Let ${C^1}, {C^1}'$ be the cubic defined by $f_3(x,y)\pm y^3=0$.
If ${C^1}, {C^1}'$ are non-singular at $P$, then $P\in C$
is an intersection singularity,  and $(C,P)$ is isomorphic to 
  $A_5, A_{11}$ or $A_{17}$
depending to the multiplicity.

Assume that $P$ is an singular point of ${C^1}$ and ${C^1}'$. 
Then $(C,P)$ can not be $E_6$  as $(C,P)$ is not irreducible
and the assertion follows from 
the classification of \cite{Pho} and   Theorem 2.\qed

Assume that $C$ is a sextics with $3A_5$ or $A_{11}+A_5$ or $A_{17}$.
We denote the location of these singularities by $P_1,P_2,P_3$
which we assume to be mutually distinct.
We say that $C$ is {\em of linear type} if
there is a line $L\subset \bfP^2$ such that 
\[
L\cap C=\begin{cases}
\{P_1,P_2,P_3\},\, I(C,L;P_i)=2\quad  & (C,P_i)\cong A_5\\
\{P_1,P_2\},\,I(C,L; P_1)=2,\, I(C,L;P_2)=4\quad &(C,P_1)\cong A_5,\,
(C.P_2)\cong A_{11}\\
\{P_1\},\,I(C,L;P_1)=6,\quad &(C,P_1)\cong A_{17}
\end{cases}
\]
The following is the converse of Proposition 7.
\begin{Proposition}\label{linear-torus}
Assume that $C$ is a reduced sextics with  $3A_5$ 
or $A_{11}+A_5$ or $A_{17}$.
Assume that 

1.  $C$ is a sextics of linear type or 

2. $C$ is a sextics of torus type which is a union of two cubics.

Then $C$ is of linear torus type.
\end{Proposition}
We give an computational proof in  Appendix.
\begin{Remark}
There exists sextics of non-torus type  which has the decomposition
type $C=B_3+B_3'$ 
with  $3A_5$ or $A_{11}+A_5$ or $A_{17}$ which are not
colinear. In fact, in the space of sextics, the moduli of sextics with
configuration $[3A_5]$, $[A_{11},A_5]$ or $[A_{17}]$ consists of  4
 components:
irreducible non-torus sextics, irreducible sextics of torus type,
   non-torus sextics with
two cubics components, sextics of linear torus type.
The assertion is shown by Artal \cite{Artal} for the configuration $[A_{17}]$.
\end{Remark}
  \section{Proof of Theorem \ref{Reduced-Simple}}
  \subsection{Elimination of other configurations}

Main step to the proof is to list the possible configurations,
 eliminating  other  configurations.
  This process can be done by fixing the inner configuration.
The proof of the existence for the survived configurations
for the maximal configurations  is  checked by constructing explicit
examples (in next subsection),
 and for other configurations, we leave it to the reader.
In the following, $B_i,B_i',\dots$ are assumed to be an irreducible
 component
of degree $i$. By \cite{Pho}, the possible inner configurations are
the combinations of $A_2,A_5,A_8,A_{11},A_{14},A_{17},E_6$.

  First consider the case $\Si_{in}(C)=[6A_2]$. 
  This implies $\de^*(C)\ge 6$. Assume that $C$ is not irreducible.
As $A_2$ is an  irreducible singularity, it is not possible unless
  $C=B_5+B_1$.
   However  there is no quintic $B_5$ with 6 $A_2$, as $n^*(B_5)=2$. 
(The conics are self-dual.)

  The configurations $\Si_{in}=[4A_2,E_6],\,
  [2A_2,2E_6],\, [3A_2,A_8],\, [A_2,A_{14}]$ are impossible to be on
 a reducible sextic curve as $\de^*\ge 7$. Now we consider the other cases.

  {\rm 1.} Assume that $\Si_{in}(C)=[4A_2,A_5]$. Then $\de^*(C)\ge 5$ and 
  the only possibility is $C=B_1+B_5$.  If this is the case,
  $B_1$ must be a flex tangent of $B_5$ and $\Si(B_5)= [4A_2]$
generically. 

  Note also  $B_1\cap B_5=A_5+2A_1$ or $A_5+A_3$ if the intersections
  are on their smooth points. 
  Here we mean by $B_1\cap B_5=A_5+2A_1$,
 the intersection of $B_1$ and $B_5$ are
  three points, and the  equivalence classes of the 
intersection  singularities are  $A_5$ and two $A_1$ respectively.
Under the assumption  $\Si_{in}(C)=[4A_2,A_5]$,
 $B_5$ can take further
$ A_1,\, A_2,\, A_3,\, A_4$, $ 2A_1$ by the genus formula.
  There are no quintic with   $4A_2+A_4$ or $5A_2+A_1$.
In fact, if there is such a quintic,
  $n^*(B_5)=3$ in both  cases. 
However this is impossible by the following well-known fact.

  Fact1.
{\em  The dual of an  irreducible smooth (resp. nodal, or cuspidal)  cubic $B_3$    is a 
9 cuspidal sextic (respectively 3 cuspidal quartic or cuspidal cubic)}.

Note that  $A_5$ mus be on $B_1\cap B_5$.
Assume first $B_1\cap B_5$ is  $A_5+2A_1$.
The configurations corresponding to the degeneration
of the quintic is:
$[A_5,4A_2,2A_1]$, $[A_5,4A_2,3A_1]$,$[A_5,5A_2,2A_1]_2$, $ [A_5,4A_2,A_3,2A_1]_1$
and $A_5,4A_2,4A_1]$.

  Assume that $B_1\cap B_5=A_5+A_3$. Then we can insert to $B_5$ either
 $A_1$ or $ 2A_1$ but we can  not insert any other singularity.
Thus we have $[A_5,4A_2,A_3]_2$,
$[A_5,4A_2,A_3,A_1]_2$ and  $[A_5,4A_2,A_3,2A_1]_2$. 
  In fact, assume that $\Si(B_5)=[5A_2]$. Then $n^*(B_5)=5$ and $i(B_5)=5$.
  Thus 5 cuspidal quintics are self dual. However
  if  $B_1\cap B_5=A_5+A_3$, $B_1$ is a flex tangent which is also
  tangent
at another point. This implies one $A_2$ of $B_5^*$ has  to be replaced
  by  $D_5$ (= the dual singularity of $( B_5, B_1\cap B_5)^*$)
  which is impossible.

The exceptional cases $[4A_2,A_5,D_4]$ and $[4A_2,A_5,D_5]$ are given
  when $\Si(B_5)=4A_2+A_1$ or $5A_2$ respectively and the line component $B_1$
  passes through the last outer $A_1$ or $A_2$. 
Note that the sextics with one of the above  configurations
can be degenerated into one of $[A_5,4A_2,A_3,2A_1]_1$ or
 $[A_5,4A_2,A_3,2A_1]_2$ or $[A_5,4A_2,D_5]$. 

There are further degenerations
$[A_5,4A_2,A_3,2A_1]_1\to [A_5,C_{3,7},A_2,A_1]$ ( 5.3-1),
$[A_5,4A_2,A_3,2A_1]_3$ 
$\to [A_5,E_6,A_3,2A_2,A_1]_2^{mr}$ (5.3-2) and 
$[4A_2,A_5,D_5]_2\to [2A_5,2A_2,D_5]_1^{mr}$ ( 5.3-3)
. We will give later explicit 
examples of these degenerations in 5.3.
So the existence  of the above configurations follows from
 the existence of these three configurations
$[A_5,4A_2,A_3,2A_1]_1$, $[A_5,4A_2,A_3,2A_1]_3$ and $[A_5,4A_2,D_5]_2$.

Note that  $[A_5,4A_2,A_3,2A_1]_i,~i=1,2$
 is an interesting weak Zariski configuration: Both has
the same decomposition
 type
$B_1+B_5$ but
 $\Si(B_5)=[4A_2,A_3]$ and $B_1+B_5=A_5+2A_1$ (respectively
 $\Si(B_5)=[4A_2,2A_1]$ and $B_1+B_5=A_5+A_3$)
for $[A_5,4A_2,A_3,2A_1]_1$ (resp. for 
$[A_5,4A_2,A_3,2A_1]_2$).
(To distinguish them, we put the index 1 or 2.)
The configurations
$[A_5,4A_2,2A_1]$, $[A_5,4A_2,3A_1]$ are also  weak
Zariski configurations as there exist irreducible sextics with 
these configurations. See \cite{Oka-Pho2}.
Hereafter we do not list up the weak Zariski configurations. They can be
 read from the indices.

  {\rm 2.} Now we consider the case 
  $\Si_{in}(C)=[2A_2,2A_5]$. \nl
(a) Consider the component type $C=B_1+B_5$.
  Then $\Si(B_5)=[2A_2,A_5]$ and $B_1\cap B_5=A_5+2A_1$ or $A_5+A_3$.
  We can put at most one $A_1$ or $A_2$ in $B_5$.
  In the case $B_1\cap B_5=A_5+A_3$, we assert that $A_2$ 
can not be inserted in $B_5$.
  In fact, assume that  $B_5$ is a quintic with
  $3A_2+A_5$.
   Then $n^*(B_5)=5$ and the dual  curve $B_5^*$
   has the same singularities,
as $i(C)=3$ and $A_5$ is self-dual
 (\cite{Oka-sextics}). If $B_1\cap B_5=A_5+A_3$, 
 $B_1$ is a flex tangent and the dual singularity is
$D_5$ which is impossible as
the dual curve $B_5^*$ can not have $A_5+2A_2+D_5$.
However the configurations $[A_5,2A_2,D_4]$ and $[A_5,2A_2,D_5]_1^{mr}$
(see 5.1-3) are 
possible by putting  the above extra $A_1$ or $A_2$ on $B_1\cap B_5$.
Note that we have a degeneration $[2A_5,3A_2,2A_1]\to [2A_5,2A_2,D_5]_1^{mr}$
\nl
(b)  Assume that $C=B_4+B_1+B_1'$. Then to have $2A_5+2A_2$,
   $B_4$ must have two cusps and  $B_1,B_1'$ must be flex tangents.
Thus generically $[2A_5,2A_2,3A_1]$. 
The configuration $[2A_5,2A_2,D_4]$ is 
given when two lines $B_1,B_1'$
intersect on $B_4$. 
Furthermore $B_4$ can have one more
node ( so $[2A_5,2A_2,4A_1]$, see 5.3-4) but it can not have three  cusp. In fact,
if $B_4$ has three cusps, $B_4^*$ is a
  nodal cubic. This is impossible as $B_4^*$ have at least two cusps.
\nl
(c)  Assume that $C=B_2+B_4$, $\Si(B_4)=[2A_2]$ and $B_2\cap B_4=2A_5+2A_1$
generically. We can put either $A_1$ or $A_2$ on $B_4$. See 5.3-5.
Consider the case $B_2\cap B_4=2A_5+A_3$.
  Then we can only insert $A_1$ into $B_4$. We can  put $A_2$ into $B_4$ 
only on $B_2\cap B_4$ so that we get $[2A_5,2A_2,D_5]_2^{mr}$ (see 5.1-3).
The case 
  $\Si(B_4)=[3A_2]$  and $B_2\cap B_4=2A_5+A_3$
does not occur. 
In fact, assume that $B_2\cap B_4=2A_5+A_3$ and $\Si(B_4)=[3A_2]$.
Note that the dual $B_2^*$ is a conic and
the dual  $B_4^*$ is cubic. Now the assumption
implies that $B_2^*\cap B_4^*=2A_5+A_3$ which is
  impossible by 
  Bezout theorem.
\nl
(d) Assume that $C=B_3+B_3'$. Then the cubics are cuspidal and
 $B_3\cap B_3'=2A_5+3A_1$ (generic) or  $2A_5+A_3+A_1$  or $3A_5$.
This is the most difficult case to find explicit examples.
  See the next section for explicit 
examples (see 5.1-1). 
The case
 $[[2A_5,2A_2],[A_5]]$
coincides with $[[3A_5],[2A_2]]$.
This corresponds to the fact that this 
configuration has two torus expressions (see 5.1-6).  Note that 
every configurations with $\Si_{in}=[2A_5,2A_2]$ except $[2A_5,2A_2,4A_1]$ is a weak Zariski
 configuration.
For example, $[2A_5,2A_2,3A_1]$ has 4 different cases.

  {\rm 3.} Assume that $\Si_{in}=[E_6,A_5,2A_2]$.
  Then $\de^*\ge 5$ and the possibility is $C=B_1+B_5$,
  $\Si(B_5)=[E_6,2A_2]$
and  $B_1\cap  B_5=A_5+2A_1$
  or $A_5+A_3$.
  As there is no quintic with $[E_6,3A_2]$ by the dual curve discussion,
  we can put at most one $A_1$. This gives the configurations
in the list. 
The  added $A_1$ can not be on $B_1$.
This has to be checked by a direct computation
or
it also follows from Yang,  \cite{Yang}, as $[E_6,A_5,2A_2,D_4]$ does
not
exist.
Note that $[E_6,A_5,2A_2,2A_1]$ and $[E_6,A_5,2A_2,A_3]$
are weak Zariski configurations.

  {\rm 4.} Assume that $\Si_{in}=[3A_5]$.
  In the case of (a)$\sim$ (c) of No.4 in Theorem \ref{Reduced-Simple},
  $B_5$ in (a), $B_4$ in (b) or (c) are rational. Thus the
  assertion
  is obvious except the existence.
The case (a) is given by $\Si(B_5)=[2A_5]$ and $B_1\cap B_5=A_5+2A_1$ 
or $[A_5+A_3]$. See 5.3-6.
The case (b) is given by $C=B_2+B_4$, $\Si(B_4)=[A_5]$ and $B_2\cap
  B_4=2A_5+2A_1$
or $2A_5+A_3$. See 5.3-7.
The case (d),  three $A_5$ are colinear by Proposition \ref{linear-torus} and assuming 
they are on $y=0$, the generic form
is given by $f_3(x,y)^2-y^6$, where $f_2=-y^2$. Thus every
  configurations
in (d) can be obtained by putting either $A_1$ or $A_2$ in the cubics.
 The cases (e,f) are special cases
of (d). In case (e), we can put only $A_1$ or $A_2$ in $B_3$.
However we need to show that if $B_1\cap
  B_2=A_3$,
   $B_3$ can not be cuspidal. In fact, 
if such a sextics exists, it gives rank 20 configuration
 $[3A_5,A_3,A_2]$
  which is known to be impossible (\cite{Shioda, Horikawa}).
  The assertion of (f) is also easy to see as three line components are 
flex tangents  and  a  nodal (respectively cuspidal)
 cubic has  three  flex points(resp. one flex point). The last
  configuration
$[3A_5,D_4]$ is given when three line components intersect at a point.

The case $(g)$ is the only non-trivial case. 
By Proposition \ref{linear-torus},
three $A_5$ can not be colinear.
The normal form is given in 5.1. In this family
 the intersection of any two conic components gives
$A_5+A_1$.
The maximal configuration $[3A_5,D_5]$ is given by $u=-1/4$
where three conics intersect at a point.

  {\rm 5.} Assume that $\Si_{in}=[2A_5,E_6]$.
  Assume that $C=B_1+B_5$ and $\Si(B_5)=[A_5,E_6]$ and therefore $B_5$ is rational.
  In the case $C=B_2+B_4$ or $C=B_1+B_1'+B_4$, $\Si(B4)=[E_6]$. In any
  case, $B_4$ is rational and the configurations
have maximal ranks. Thus there are no further possibility. See 5.1-4,5.

  {\rm 6.} Assume that $\Si_{in}(C)=[A_8,A_5,A_2]$. Then  $C=B_1+B_5$.
The non-existence of $\Si(B_5)=A_8+2A_2$ with
$\Si_{red}=[A_8,A_5,2A_2,2A_1]$ 
is checked by a direct computation. This is also a result of Yang
\cite{Yang},
as this is not in his table of maximal rank configuration.
The  $[A_8,A_5,2A_2,A_3]$
 does not exists as  the rank is 20.

  {\rm 7.} Assume that $\Si_{in}(C)=[A_{11},2A_2]$.
It is easy to see that $A_{11}$ must be an intersection singularity.
\nl
(a) Assume first $C=B_2+B_4$.
As $B_4$ has $2A_2$, it can take only $A_1$ or $A_2$.
The case $B_2\cap B_4=A_{11}+A_3$, it can be checked by computation that 
$B_4$ can not have further singularity.  Namely $[A_{11},A_3,2A_2,A_1]$ 
does not appear from this series.
 This  also follows from  the connectedness of the moduli
space of the sextics with  the configuration  $[2A_2,A_{11},A_3,A_1]^{mr}$
  (see \cite{Yang}), as it exists  for the component type
$C=B_3+B_3'$. See 5.3-8.
\nl
(b)
Assume that  $C=B_3+B_3'$. Then two cubics are cuspidal
and $B_3\cap B_3'=A_{11}+3A_1$
or $A_{11}+A_3+A_1$. 
As $[A_{11},2A_2,A_3,A_1]^{mr}$ (see 5.1-10) has rank 19, there
  are no  possibility of $[A_{11},A_5,2A_2]$.


{\rm 8.} Assume that $\Si_{in}=[A_{11},A_5]$. 
\nl
(a) Assume that 
 $C=B_1+B_5$ and
$\Si(B_5)=[A_{11}]$. Then $B_4$ is rational and
 we can not put any  further singularity in $B_5$.
As $B_1\cap B_5=A_5+2A_1$ or $A_5+A_3$, the assertion is clear.

Assume that $A_{11}$ is an intersection singularity. Then it implies
either
$C=B_2+B_4$, $C=B_3+B_3'$ or $C=B_1+B_2+B_3$.
\nl
(b) Assume that $C=B_2+B_4$. Then $\Si(B_4)=[A_5]$ and $B_2\cap
B_4=A_{11}+2A_1$ or $A_{11}+A_3$ and the assertion is clear.
\nl
(c) Assume that $C=B_3+B_3'$. We can make two cubics are tangent at two
points
with multiplicity 6 and 3 so that $B_3\cap B_3'=A_{11}+A_5$.
Now the assertion follows by putting $A_1$ or $A_2$ in the cubics.
As the rank is 19 for $[A_{11},A_5,A_2,A_1]$, we can not put $A_2$ in
both cubics simultaneously.
\nl
(d)
Assume that $C=B_3+B_2+B_1$. We can make them so that 
$B_3\cap B_2=A_{11}$, $B_1\cap B_3=A_5$ and $B_1\cap B_2$ is either 
$2A_1$ or $A_3$ and $B_3$ is either smooth or nodal.

{\rm 9.} Assume that $\Si_{in}=[A_{17}]$. Then th only possibility is 
$C=B_3+B_3'$ with $B_3\cap B_3'=A_{17}$.
 The assertion is clear. See 5.1-17,18.

\section{Proof of Theorem \ref{Reduced-NonSimple}}
In this section, we prove Theorem \ref{Reduced-NonSimple}.
As in the proof of Theorem \ref{Reduced-Simple}, the proof of existence
after eliminating non-existing configurations, is due to direct
computations
and we give some non-trivial examples later.
We first fix a non-simple inner singularity at the origin $O$
and then we consider the possibility of inner configurations and 
component types.

\noindent
1. Assume that $B_{3,6}\in \Si_{in}$. Recall that $(C_2,O)$ is smooth,
$(C_3,O)\cong A_1$ and $\iota:=I(C_2,C_3;O)=3$.
Possible inner configurations are 
$ [3A_2,B_{3,6}],\, [A_2,A_5,B_{3,6}]$,
$[A_2,E_6,B_{3,6}]$, $[A_8,B_{3,6}]$, $[2B_{3,6}]$.
We assume that  $B_{3,6}$  is at $O$.
First observe that has locally
three smooth components $C^1,C^2,C^3$ with $I(C^i,C^j;O)=2$ for $i\ne j$. Thus 
  if $B_{3,6}$ is an intersection singularity
of two  global components, say $C^1, C^2\cup C^2$, $(C^2\cup C^3;O)\cong A_3$
 and  $I(C^1,C^2\cup C^3;O)=4$.
\nl
(a) Assume that  $\Si_{in}=[3A_2,B_{3,6}]$. If $C$has two components,
the singularities $3A_2, A_3$ must be in a component. Then the unique
      possibility is  the case $C=B_1+B_5$ and $B_{3,6}$ is an intersection
      singularity of $B_1$ and $B_5$ and 
 $B_5$ has 
      $A_3+3A_2$  as singularities. 
Thus    $I(B_1,B_5;O)=4$, $B_1\cap B_5=B_{3,6}+A_1$ and
 $B_5$  can take further at most
      either one
$A_1$ or   $A_2$. It is easy to observe that $C$ can not have three
irreducible components  (see 5.2-1).
\nl
(b)
Assume that $\Si_{in}=[B_{3,6},E_6,A_2]$ or $[B_{3,6},A_8]$. Then 
by an easy consideration about $\de^*$-genus, $C=B_1+B_5$ and $\Si_{B_5}$
is $[A_2,E_6,A_3]$ or $[A_8,A_3]$ and 
 $B_5$is already rational. We get
 $\Si_{red}=[B_{3,6},E_6,A_2,A_1],\,[B_{3,6},A_8,A_1]$.
 See 5.2-4.
\nl
(c)
Assume that $\Si_{in}=[A_2,A_5,B_{3,6}]$.  In this case,
we have more possibilities of component types.
\nl
(1) Assume that $C=B_1+B_5$. Then
$\Si({B_5})= [A_2,A_5,A_3]$  and   $B_5$ is
      rational.
\nl
(2) Assume that $C=B_2+B_4$. Then $\Si(B_4)=A_3+A_2$ and $B_4$ is thus
  rational and $B_2\cap B_4= B_{3,6}+A_5+A_1$.
\nl
(3) Assume that $C=B_1+B_1'+B_4$. $B_4$ is as above and $B_1\cap B_4=B_{3,6}$
and $B_1'\cap B_4=A_5+A_1$. See 5.2-2.
\nl
(4) Assume that $C$ has  a cubic component $B_3$. Then $B_3$ has to take
$A_2$ and smooth at $O$.
To make $B_{3,6}$, the other components can not be three lines, 
As an
      irreducible
cubic can not have $A_3$, the only possibility is $C=B_1+B_2+B_3$,
$B_2\cap B_3=A_3+A_5+A_1$ and $B_1$ is tangent to $B_2$ and $B_3$ at 
$O$ so that altogether we  get $B_{3,6}$ and $B_1\cap B_3$ 
has one transverse intersection to make $A_1$. See 5.2-3.
\nl
(5) Finally for  $\Si_{in}=[2B_{3,6}]$, it is already observed in \cite{Pho} that
      $C=B_2+B_2'+B_2''$.

\vspace{.3cm}
\noindent
2. Assume that $C_{3,7}\in \Si_{in}$. By the classification \cite{Pho},
 $C_2$
is smooth, $C_3$ is nodal at $O$ and $\iota=3$.
 Recall that $C_{3,7}$ is an
      intersection
singularity of a smooth component and a component with $A_4$.
This implies that $C$ must have a component of degree $\ge 4$.
The possibilities of $\Si_{in}$ are
$ [3A_2,C_{3,7}]$, $ [A_2,A_5,C_{3,7}]$, $ [A_2,E_6,C_{3,7}]$, 
$[A_8, C_{3,7}]$
(\cite{Pho}). In any case, as $\de(C_{3,7})=6$, $C_{3,7}$ must be  an
intersection singularity. We assume that $O$ is $C_{3,7}$-singularity.
\nl
(a) Assume that $\Si_{in}= [3A_2,C_{3,7}]$. 
Then $C=B_1+B_5$. Note that $I(B_1,B_5;O)=4$ by Example \ref{example-inter} in 
 the section 2 and  $B_1\cap B_5=C_{3,7}+A_1$.
Then 
$ \Si(B_5)\supset 3A_2+A_4$, we can put $A_1$ in $B_5$.
Note that a quintic $B_5$ can not have  $4A_2+A_4$ as, if so, 
we get   $n^*(B_5)=3$ which is a contradiction.
\nl
(b) Assume that $\Si_{in}=  [A_2,E_6,C_{3,7}]$ or $[A_8, C_{3,7}]$.
Then the possibility is $C=B_1+B_5$ and $\Si(B_5)=A_2+E_6+A_4$
or $A_2+E_6+A_4$. In any case $B_5$ is rational and $B_1\cap
B_5=C_{3,7}+A_1$ and we get $[C_{3,7},E_6,A_2,A_1]$ and
$[C_{3,7},A_8,A_1]$.
See 5.2-5.
\nl
(c)Assume that  $\Si_{in}=[A_2,A_5,C_{3,7}]$.
\nl (1) If $C=B_1+B_5$,
 $\Si(B_5)=A_2+A_5+A_4$ and thus $B_5$ is rational.
Thus $\Si_{red}=[C_{3,7},A_5,A_2,A_1]$.
\nl
(2) Assume that $C=B_2+B_4$.
 Then $\Si(B_4)=A_4+A_2$ and thus  $B_4$ is rational
and $B_2\cap B_4=A_5+C_{3,7}+A_1$ as $I(B_2,B_4;O)=4$.
Assume that $C=B_1+B_1'+B_4$. Then $B_4$ is as above and $B_1\cap B_4=C_{3,7}$
and $B_1'\cap B_4=A_5+A_1$ and the corresponding configuration is 
$[C_{3,7},A_5,A_2,2A_1]$. See 5.2-6.

\par\vspace{.3cm}\noindent
3.  Assume that $C_{3,8}\in \Si_{in}$. Then  $C_2$ is smooth,
$C_3$ is nodal at $O$ and $\iota=3$.
Assume that $O$ is $C_{3,8}$ singularity defined by 
$y^3+y^2x^2-x^8+\text{(higher terms)}=0$ for simplicity.
 Recall that it has three smooth 
components $L_1,L_2,L_3$
where $L_1: y+x^2+\text{(higher terms)}=0$
and $L_2,L_3: y\pm x^3+\text{(higher terms)}=0$. To consider it as an intersection singularity
of two components,
there are two ways.
\nl
(a-1) Assume that $L_2$ is a smooth component of $C$ and $L_1\cup L_3$ is
      another component. Then $I(L_2,L_1\cup L_3; O)=5$ and 
$(L_1\cup L_3;O)=A_3$.
\nl
(a-2)  Assume that $L_1$ is a smooth component of $C$ and $L_2\cup L_3$ is
      another component. Then $I(L_1,L_2\cup L_3; O)=4$ and 
$(L_2\cup L_3;O)=A_5$.

Possible inner configurations are
$[3A_2,C_{3,8}]$, $ [A_2,A_5,C_{3,8}]$, $[A_2,E_6,C_{3,8}]$, $
[A_8,C_{3,8}]$.
\nl (1) Assume that $\Si_{in}=[3A_2,C_{3,8}]$ or 
 $[A_2,E_6,C_{3,8}]$ or $[A_8,C_{3,8}]$. Then $C=B_1+B_5$.
\nl
1. Assume first that   $B_1$ corresponds to $L_2$ and $B_5$
      corresponds to
$L_1\cup L_3$.
Then $(B_5,O)\cong A_3$, $B_1\cap B_5=C_{3,8}$ and generically we have
$ \Si({B_5})=\, [ 3A_2,A_3]$, $ [A_2,E_6,A_3]$,$ [A_8,A_3]$ respectively.
In the last two cases, $B_5$ is rational and
it is easy to see that (a-2) does not occur. 
We get $\Si_{red}=[C_{3,8},E_6,A_2]$ and $[C_{3,8},A_8]$.

Assume that $\Si_{in}=[3A_2,C_{3,8}]$. Then $\Si(B_5)=[3A_2,A_3]$
and we can put further $A_1$. Thus we get 
$\Si_{red}=[C_{3,8},3A_2],\,[C_{3,8},3A_2,A_1] $. See 5.2-7.
We assert that we can not put $A_2$ in $B_5$:
\begin{Assertion} Such a quintic $B_5$ with $\Si(B_5)=[4A_2,A_3]$  does not exist.
\end{Assertion}
{\em Proof.}
Suppose that such a quintic exists. Then $n^*(B_5)=4$.
By the assumption, $(C,O)$ has locally three components $L_1,L_2, L_3$
$B_5$ has  locally two components $L_1, L_3$.
 Recall that $I(L_2,L_3;O)=3$.
As we have assumed that $B_1=L_2$, this implies that $L_3$ has a flex
      point
at $O$. Other component $L_1$ has $I(L_2,L_1;O)=2$.
 Assuming $(x,y)$ is an affine
coordinate system so that $y=0$ be the equation of $L_2$, $L_1$ and $L_2$ are defined by
$h_1(x,y)=(y+a\,x^2+\text{(higher terms)})$ and
$h_3(x,y)=(y+\,b\, x^3+\text{(higher terms)})=0$
for some $a,\, b\ne 0$. Here $h_1,h_2$ are analytic functions defined 
in a neighborhood of $O$, though $(x,y)$ are  affine coordinates.
This implies by the following lemma that  the dual singularity of  $(B_5,O)$
is a union of a cusp $L_3^*$ and a smooth curve $L_1^*$
which has the same tangent 
with the cusp. Thus   the Milnor number of $(B_5^*,O^*)$ is 7.
(This implies $A_3$ is not generic in the sense of Puiseux order
     \cite{Oka-sextics}).
 However a quartic can 
have at most 6 as the total Milnor number, which is a contradiction.\qed
\begin{Lemma}
Let $B_5$ be a projective curve with a singularity at the origin
whose defining function takes the form
$h_1(x,y)h_3(x,y)$.
Then the dual singularity $(B^*,O^*)$ is locally defined 
by $g(u,v)=0$ where 
\begin{eqnarray*}
& g(u,v)=(v+a'\, u^2+\text{(higher terms)})(v^2+b'\, u^3+\text{(higher terms)}),
\quad a',\, b'\ne 0\\
&=v^3+b'\, v u^3\,+\, a'\,b' u^5+\text{(higher term)}=0\qquad\qquad\qquad\qquad
\end{eqnarray*}
Thus  the dual singularity is $E_7$ and
the Milnor number  is 7.
\end{Lemma}
{\em Proof.} We use the parametrization
$L_1:\,x(t)=t,\, y(t)=-at^2+\text{(higher terms))}$
and $L_3:\, x(t)=t,\, y(t)=-b\, t^3+\text{(higher terms)}$.
Then the equation of the Gauss map images can be obtained by 
an easy computation (see \cite{Oka-sextics}) and the assertion follows.\qed

2. Now we consider the case $C=B_1+B_5$ which corresponds to (a-2):
$(B_5,O)\cong A_5$ and $B_1\cap B_5=C_{3,8}+A_1$. Then
the unique possible inner configuration is 
$\Si_{in}=[C_{3,8},3A_2]$ and $\Si(B_5)=[3A_2,A_5]$ and $B_5$ is
rational. 
This gives the
configuration $[C_{3,8},3A_2,A_1]_2$. See 5.2-8.
\nl
(2) Now assume that $\Si_{in}=[ C_{3,8},A_5,A_2]$.
\nl
1. Assume that $C=B_1+B_5$. As $\Si(B_5)\supset \{A_5,A_2\}$, $B_5$ can
not take another $A_5$. Thus  the case (a-1) is  the unique
possibility
and $B_5$ is rational with $[A_5,A_3,A_2]$. This gives the configuration
$\Si_{red}=[C_{3,8},A_5,A_2]$.
\nl
2. Assume
 $C=B_2+B_4$. As we have seen in (a-1) and (a-2), we need
either $A_3$ or $A_5$ on $B_4$ to make $C_{3,8}$. Thus the only
      possibility
is the case $\Si(B_4)=[A_3,A_2]$  ((a-1)) and 
$\Si_{in}=[A_2,A_5,C_{3,8}]$. In fact, this case is possible  and 
$B_2\cap B_4=C_{3,8}+A_5$ as $I(B_2,B_4;O)=5$
and $\Si_{red}=[A_2,A_5,C_{3,8}]$.
\nl
3. Now we consider the case $C$ does not have any component of degree
      greater than 3. 
Only possible inner configuration is $[A_2,A_5,C_{3,8}]$
and the component type must be $\{1,2,3\}$ and the cubic
      component
must be cuspidal.
To make $C_{3,8}$, we need either $A_3$ or $A_5$ in the other union of 
components.  We can make the components $B_1,B_2,B_3$ in two ways.
\nl
(a-1) $I(B_1,B_2;O)=2$, $I(B_1,B_3;O)=3$ and $I(B_2,B_3;O)=2$ and 
$(B_2\cup B_3,O)\cong A_3$: The corresponding configuration is denoted
by $[C_{3,8},A_5,A_2,A_1]_1$. This is a degeneration of $[C_{3,8},A_5,A_2]_1$.
See 5.2-9.
\nl
(a-2) $I(B_1,B_2;O)=2$, $I(B_1,B_3;O)=2$ and $I(B_2,B_3;O)=3$ and 
$(B_2\cup B_3,O)\cong A_5$: The corresponding configuration is denoted
by $[C_{3,8},A_5,A_2,A_1]_2$. This is a degeneration of $[C_{3,8},A_5,A_2]_2$.
See 5.2-10.
\par\vspace{.3cm}\noindent
4. $C_{3,9}\in \Si_{in}$. Then $C_2$ is smooth and $C_3$ is nodal at $O$
and $\iota=3$ or $4$. We assume that $O$ is 
$C_{3,9}$-singularity as before.
 First we observe that $\mu(C_{3,9})=13$ and it must be
 an intersection singularity of a smooth component $L$
and a component  $M$ with $A_6$. Note that $I(L,M;O)=4$.
There are two $C_{3,9}$ with different
      $\iota$ (=the intersection number  $I(C_2,C_3;O)$).
\nl 1.
The case $\iota=3$, the only possibility of the inner configuration is 
 $\Si_{in}=[3A_2,C_{3,9}]$ by \cite{Pho} which is given by  $C=B_1+B_5$ and 
$\Si(B_5)=3A_2+A_6$ and $B_1\cap B_5=C_{3,9}+A_1$,  and
      $\Si_{red}=[3A_2,C_{3,9},A_1]$.
\nl 2.
The other case is $\iota=4$ and the possible inner configurations are 
$ [2A_2,C_{3,9}]$, $ [A_5,C_{3,9}]$, $ [E_6,C_{3,9}]$.
\nl
a. Assume $\Si_{in}=[C_{3,9},2A_2]$. Then $\Si(B_5)=2A_2+A_6$ and  we can put
$A_1$ or $A_2$. 
This gives $\Si_{red}=[C_{3,9},2A_2,A_1],\, [C_{3,9},2A_2,2A_1], \,
[C_{3,9},3A_2,A_1]$.
\nl b.
In the cases $\Si_{in}=[A_5,C_{3,9}]$ or $ [E_6,C_{3,9}]$,
we have  $\Si(B_5)=A_5+A_6$ or $E_6+A_6$
and therefore $B_5$ is already rational.
Note that $I(B_1,B_5; O)=4$ and thus $B_1\cap B_5=C_{3,9}+A_1$.
Thus the possibilities for $\Si_{red}$ are
$[2A_2,C_{3,9},A_1]$, $ [2A_2,C_{3,9},2A_1]$, 
$[3A_2,C_{3,9},A_1]$, $ [A_5,C_{3,9},A_1]$, $ [E_6,C_{3,9},A_1]$.
See 5.2-12.
\begin{Remark}
We get the reduced configuration $[C_{39},3A_2,A_1]$
from two inner configurations $[C_{3,9}, 3A_2]$ and
$[C_{3,9},2A_2]$
In fact, the moduli is the same and it has two different
torus decompositions. An example is the following.
\begin{multline*}
\mathrm{f}(x, \,y) := (y^{2} + (x + 1)\,y - x^{2})^{3} + (y^{3}
 + ({\displaystyle \frac {16}{3}} \,x + 1)\,y^{2} + (6\,x^{2} + 3
\,x)\,y + x^{3})^{2}
\\
={\displaystyle \frac {343}{27}} \,(y^{2} + {\displaystyle \frac {
15}{7}} \,y\,x + {\displaystyle \frac {3}{7}} \,y + 
{\displaystyle \frac {9}{7}} \,x^{2})^{3} - {\displaystyle 
\frac {1}{27}} \,(27\,x^{3} + 9\,y\,x + 60\,y\,x^{2} + 9\,y^{2}
 + 54\,y^{2}\,x + 17\,y^{3})^{2}
\end{multline*}
\end{Remark}

\par\vspace{.3cm}\noindent
5. $C_{3,12}\in \Si_{in}$. In this case, $C_2$ is smooth, $C_3$ is nodal
at $O$ and $\iota=5$. Possible
inner configuration is $[A_2,C_{3,12}]$. 
First note that $C_{3,12}$ has locally three smooth components
$L_1,L_2,L_3$ which satisfies
\begin{eqnarray*}
& I(L_1,L_2;O)=I(L_1,L_3;O)=2, \, I(L_2,L_3;O)=5\\
&(L_1\cup L_2;O),\, (L_1\cup L_3;O)\cong A_3,\,\quad (L_2\cup L_3;O)\cong A_9
\end{eqnarray*}
If $C$ has two components, it can be either
$L_1+(L_2\cup L_3)$ or $L_2+(L_1\cup L_3)$.

Assume that $C=B_1+B_5$. 
As $I(B_1,B_5;O)\le 5$, we must have  $B_5=L_2\cup L_3$. Then
$B_1\cap B_5=C_{3,12}+A_1$ and $\Si(B_5)=[A_2,A_9]$ and $B_5$ is
      rational. This gives the configuration $[C_{3,12},A_2,A_1]_1$.
See 5.2-13.

Assume that $C=B_2+B_4$. Then as $B_4$ can not have $A_9$,
we must have $B_4=L_1\cup L_3$
and  $\Si(B_4)=[A_3,A_2]$ and 
$B_2\cap B_4=C_{3,12}+A_1$ as $I(L_2,L_1\cup L_3;O)=7$.
Thus  $B_4$ is rational and
$\Si_{red}=[A_2,C_{3,12},A_1]_2$. See 5.2-13.

Assume that $C_{3,12}$ is an intersection singularity
of  three global components.
The intersection singularity of two of them
  have to make   $A_9$.
To make $A_9$, we need the intersection multiplicity 5.
Thus the unique possibility is the case:
$C=B_1+B_2+B_3$ with $(B_2\cup B_3,O)\cong A_9$. 
Thus we may assume that $B_1=L_1,\, B_2=L_2,\, B_3=L_3$,
 $B_2\cap B_3=A_9+A_1$ and 
$B_1$ is tangent to $B_2$
at $O$ so that $(B_1\cup B_2\cup B_3,O)\cong C_{3,12}$. The corresponding
configuration is $[C_{3,12},A_2,2\,A_1]$. 
This is a degeneration of  $[C_{3,12},A_2,A_1]_i,\, i=1,2$.

\par\vspace{.3cm}\noindent
6. $C_{6,6}\in \Si_{in}$. In this case, both of $C_2$ and $C_3$ are nodal
at $O$ and $\iota=4$. Possible inner configurations are
$[2A_2,C_{6,6}]$ and $[A_5,C_{6,6}]$. We assume   $C_{6,6}$
      singularity
 is at $O$ and is locally defined by
$y^6-x^2y^2+x^6+\text{(higher terms)}=0$ for simplicity.
First note that $C_{6,6}$ has locally 4 smooth components
$L_1,L_2,K_1,K_2$ such that $L_1,L_2: x\pm y^2+\text{(higher terms)}=0$
and $K_1,K_2: y\pm x^2+\text{(higher terms)}=0$ and
$I( K_1,K_2;O)=I(L_1,L_2;O)=2$ and $ I(L_i,K_j;O)=1$.
Note also that $(L_1\cup L_2\cup K_1,O)\cong D_6$
(same for any three components) and $(L_1\cup L_2, O)$, 
\nl
$ (K_1\cup K_2,O)\cong A_3$.
\nl
1. Assume that $C=B_1+B_5$. Then $\Si(B_5)=2A_2+D_6$ and
      $\Si_{red}=[2A_2,C_{6,6},A_1]$.
As $B_5$ is already rational, the case $\Si_{in}=[C_{6,6},A_5]$ does not occur.
\nl 
The case $C=B_2+B_4$ does not exist because
 $B_4$ can not support $ D_6$.
\nl
2. Assume that $C=B_3+B_3'$.
Then this case is possible only if
 $\Si_{in}=[A_5,C_{6,6}]$ and two cubics are nodal
at $O$ with the same tangent cone so that $B_3\cap B_3'=A_5+C_{6,6}$ 
and $\Si_{red}=[A_5,C_{6,6}]$.
Note that $I(L_1\cup K_1,L_2\cup K_2;O)=6$.
Another decomposition possibilities:
\nl 3.
Assume $C=B_4+B_1+B_1'$: $\Si(B_4)=[2A_2,A_1]$ and $ B_1$  and $B_1'$
are tangent to branches of $A_1$ of $B_4$ so that 
$(B_1\cup B_1' \cup B_4,O)=C_{6,6}$
and $\Si_{red}=[2A_2,C_{6,6},2A_1]$. Two $A_1$ are
the transverse intersection of $B_1$ or $B_1'$ and $B_4$ outside of $O$.
See 5.2-14.
 As $I(B_1,B_4;O)=3$, the case 
$\Si_{in}=[C_{6,6},A_5]$ does not exist.
\nl 4.
Assume $C=B_3+B_2+B_1$: $B_3$ is nodal and $\Si_{in}=[A_5,C_{6,6}]$.
 $B_1\cap  B_2\cap B_3=C_{6,6}$ and 
 $B_2\cap B_3=A_5+D_6$ and $\Si_{red}=[C_{6,6},A_5]$.
\nl 5.
$C=B_2+B_2'+B_1+B_1'$: This can be understood as a degeneration of 
$B_3+B_3'$ and $\Si_{red}=[A_5,C_{6,6},2A_1]$. See 5.2-15.
The case $C=B_3+B_1+B_1'+B_1''$ can not make $C_{6,6}$.

\par\vspace{.3cm}\noindent
7. $C_{6,9}\in \Si_{in}$. Possible inner configuration is
       $[A_2,C_{6,9}]$.
 Note that $C_{6,9}$ has two smooth components $L_1,L_2$
defined by
$L_i: y+a_i x^2+\text{(higher terms)}=0$, $a_i\ne 0,a_1\ne a_2$, and one
       component $K$ defined by  $y^2+bx^7+\text{(higher terms)}=0,b\ne 0$
 with $A_6$ singularity and $I(L_1,L_2;O)=2$ and $I(L_i,K;O)=2$.
 Note that $(L_2\cup K;O)\cong D_9$.
 Thus  the unique possibility  is the case $C=B_1+B_5$
 with $\Si(B_5)=[A_2,D_9]$, $B_1\cap B_5=C_{6,9}+A_1$ and
       $\Si_{red}=[A_2,C_{6,9},A_1]$.

 \par\vspace{.3cm}\noindent
 8. $C_{3,15}\in \Si_{in}$. Then $C_2$ is smooth and $C_3$ is nodal at $O$.
 Note that $C_{3,15}$ has two components:
 a smooth component $L$ and another component $K$ with $A_{12}$ singularity
       and
 $I(L,K;O)=4$. Thus  the unique possibility is the case $C=B_1+B_5$,
       $\Si(B_5)=[A_{12}]$
 and $B_1\cap B_5=C_{3,15}+A_1$. This case gives
       $\Si_{red}=[C_{3,15},A_1]$. See 5.2-17.

 \par\vspace{.3cm}\noindent
 9. $B_{3,12}\in \Si_{in}$: This case is unique and $C=B_2+B_2'+B_2''$
and
$\Si_{red}=[B_{3,12}]$. See 5.2-18 and \cite{Pho}.
 \par\vspace{.3cm}\noindent
10. $C_{6,12}\in \Si_{in}$. In this case,  $C_2$ is a multiple line and 
$C_3$ is 
nodal at $O$. Note that $C_{6,12}$ has 4 smooth
       components
 $L_1,L_2,K_1,K_2$ with $I(L_1,L_2;O)=2$, $I(K_1,K_2;O)=5$
 and $I(L_i,K_j;O)=1$. Note also that 
\begin{eqnarray*}
& (L_1\cup K_1;O)\cong A_1,\,(L_1\cup L_2;O)\cong A_3,\,
K_1\cup K_2;O)\cong A_9\\
& (L_1\cup L_2\cup K_1;O)\cong D_6,\,(L_2\cup K_1\cup K_2,O)=D_{12}
\end{eqnarray*}
 Thus  $C=B_2+B_4$ is not possible. If $C=B_1+B_5$ is the case, 
$(B_5,O)\cong D_6$ and
$B_5=L_1\cup L_2\cup K_1$.
But this is impossible as  $5\le I(B_1,B_5;O)=7$.

 Assume that $C=B_3+B_3'$. Then the cubics are nodal and they  correspond to 
$L_i\cup K_i, \, i=1,2$ respectively
 and 
 $I(B_3,B_3;O)=9$.  This case exists and  $\Si_{red}=[C_{6,12}]$.

 Suppose that $B_3'$ degenerate into $B_1+B_2$.
 Then $B_1\cap B_3=D_6$ and $B_2\cap B_3=D_{12}$ and $B_1\cap B_2=2A_1$.
 This also exists and  $\Si_{red}=[C_{6,12},A_1]$.
See 5.2-16.

There are no other  possibility of three components case.
 Note also that
       $C$
 can not have two line components. In fact, if it has two line
 components,
we may assume that 
       $L_1,K_1$ are  the lines components.
Put $C=L_1+K_1+J$, where $J$  is the union of other components,
$\degree(J)=4$ and  $J=L_2\cup K_2$.
Then we have a contradiction
 $8\ge I(L_1\cup K_1,J;O)= I(L_1\cup K_1, L_2\cup K_2;O)=9$.

 \par\vspace{.3cm}\noindent
 11. $B_{4,6}\in \Si_{in}$.
In this case, $C_2$ is a multiple line and $C_3$ has a cusp ( or 
$A_3$). The possible inner configurations are $[2A_2,B_{4,6}]$
 and $[A_5,B_{4,6}]$. Note that $B_{4,6}$ can be intersection singularity
       of two cuspidal components
 with intersection number 6. Thus the only possibility is that 
 $C=B_3+B_3'$ and it is easy to see that $[2A_2,B_{4,6}]$ does not exist
as a reducible sextics.
 For $[A_5,B_{4,6}]$, we can take two cuspidal cubic $B_3,B_3'$
 such that $B_3\cap B_3'=B_{4,6}+A_5$.  There are no other possibility.

 \par\vspace{.3cm}\noindent
 12. $D_{4,7}\in \Si_{in}$. Then  the possible inner configurations are
 $[2A_2,D_{4,7}]$ and $[A_5,D_{4,7}]$.
 Recall that  $D_{4,7}$ is defined by 
 $y^4+x^3y^2+a y x^5+b x^7=0$ with $a^2-4b\ne 0$ and $\mu(D_{4,7})=16$. It has three components
 $L_1,L_2,K$ where $K$ is cuspidal component of
type $x^3+y^2+\text{(higher terms)}=0$ and $L_1,L_2$ are smooth
components of type $y+\al x^2+\text{(higher terms)}=0$ and thus 
       $I(L_1,L_2;O)=2$
 and $I(L_i,K;O)=3$. Thus  $(L_1\cup L_2;O)\cong A_3$
and $(L_1\cup K,O)\cong E_7$.
 For $[2A_2,D_{4,7}]$, as an component have to support 
$2A_2$ and  $E_7$, the only possibility is the case $C=B_1+B_5$
 and $\Si(B_5)=2A_2+E_7$ and  so $B_5$ is rational and $B_1\cap
 B_5=D_{4,7}$.

Consider the case $\Si_{A_5}=[A_5,D_{4,7}]$. If $D_{4,7}$
is an intersection singularity of two components, this is only
possible for $C=B_1+B_5$, but then we can not make $A_5$.
Thus $C$ has three components and   the unique possibility is
       $C=B_3+B_2+B_1$.
 In fact, this is possible if $B_3$ is cuspidal
 and $B_2\cup B_3=E_7+A_5$ and $\Si_{red}=[A_5,D_{4,7}]$.
See 5.2-19 and  \cite{Pho}.

 \par\vspace{.3cm}\noindent
 13. $Sp_2\in \Si_{in}$. This  case is studied by \cite{Pho} and given
       by 
 $C=B_3+B_3'$, where both cubics are cuspidal with the same tangent cone
and $I(B_3,B_3';O)=9$.

 \par\vspace{.3cm}\noindent
 14. $B_{6,6}$ is possible only when $f_2(x,y), f_3(x,y)$ are
       homogeneous polynomials of degree 2 and 3 in $x,y$.

 \subsection{Examples II, Non-simple singularities}
 We give explicit examples for some configurations.
 The normal form of $B_{3,6}$-singulaity at $O$ with 
 $y=0$ a linear component is given as 
\begin{eqnarray*}
&\mathit{f_2} := \mathit{a_{02}}\,y^{2} + (\mathit{a_{11}}\,x + \mathit{
a_{01}})\,y - t^{2}\,x^{2}
\\
&\mathit{f_3} := \mathit{b_{03}}\,y^{3} + (\mathit{b12}\,x + \mathit{
b_{02}})\,y^{2} + (\mathit{b_{21}}\,x^{2} + \mathit{b_{11}}\,x)\,y + t^{3}
\,x^{3}
\end{eqnarray*}
The normal form of torus decomposition with $C=B_2+B_4$ and $B_{3,6}$ at
$O$
where $B_2$ is defined by $y-x^2=0$ is given
\begin{multline*}
\mathrm{f}(x, \,y) := ( - t^{2}\,y^{2} + (\mathit{a_{11}}\,x + 
\mathit{a_{01}})\,y - {\displaystyle \frac {1}{4}} \,{\displaystyle 
\frac {(4\,t^{2}\,\mathit{a_{01}} + \mathit{a_{11}}^{2})\,x^{2}}{t^{2}}
} )^{3} + {\displaystyle \frac {1}{64}} ( - 8\,t^{6}\,y^{3} \\
\mbox{} + 12\,y^{2}\,t^{4}\,\mathit{a_{11}}\,x - 8\,y^{2}\,t^{3}\,
\mathit{b_{02}} - 6\,x^{2}\,y\,t^{2}\,\mathit{a_{11}}^{2} + 8\,x^{2}\,y
\,t^{3}\,\mathit{b_{02}} - 8\,x\,y\,t^{3}\,\mathit{b_{11}} \\
\mbox{} + x^{3}\,\mathit{a_{11}}^{3} + 8\,x^{3}\,t^{3}\,\mathit{b_{11}}
)^{2}/t^{6} 
\end{multline*}
In the following examples, those with  a line or conic component can be easily
derived from
the above normal forms.
 \begin{enumerate}
\item $B_{3,6}\in \Si_{in}$:
$\Si_{red}=[4A_2,B_{36},A_1]$, $C=B_1+B_5,\, \Si(B_5)=[4A_2,A_3]$:

\begin{multline*}
\mathrm{f}(x,y)=
 (y^{2} + ( - 5\,x + 1)\,y - x^{2})^{3} + (
 - {\displaystyle \frac {3}{4}} \,y^{3} + ({\displaystyle \frac {
15}{4}} \,x + {\displaystyle \frac {3}{4}} )\,y^{2} + ( - 15\,x^{
2} + {\displaystyle \frac {9}{4}} \,x)\,y + x^{3})^{2}
\end{multline*}
\item  $\Si_{red}=[A_2,A_5,B_{3,6},2A_1]$, $C=B_1+B_1'+B_4$:
\[
 \mathrm{f}(x,y)
= ( - y^{2} + y - 4\,x^{2})^{3} + (y^{3} + ( - 4\,x - 1)\,y^{2} + 4\,y\,x - 8\,x^{3})^{2}
\]
\item  $\Si_{red}=[A_2,A_5,B_{3,6},2A_1]$, $C=B_3+B_2+B_1$:
\[
 \mathrm{f}(x,y)
= ( - y^{2} + ( - 3\,x + 1)\,y - x^{2})^{3}
 + ({\displaystyle \frac {2}{3}} \,y^{3} + ({\displaystyle 
\frac {13}{3}} \,x - {\displaystyle \frac {2}{3}} )\,y^{2} + (5\,
x^{2} - {\displaystyle \frac {7}{3}} \,x)\,y + x^{3})^{2}
\]
\item  $[A_8,B_{3,6},A_1]$, $C=B_1+B_5$ and $\Si(B_5)=[A_8,A_3]$:
\[
 \mathrm{f}(x,y)=
( - y^{2} + ( - 3\,x + 1)\,y - x^{2})^{3}
 + ( - y^{3} + ({\displaystyle \frac {2}{3}} \,x + 1)\,y^{2} + (
10\,x^{2} - {\displaystyle \frac {11}{3}} \,x)\,y + x^{3})^{2}
\]
\item $\Si_{red}=[C_{3,7},A_8,A_1]$ with $C=B_1+B_5$:
\[
\mathrm{f}(x, \,y) := ( - y^{2} + y - x^{2})^{3} + ( - 2\,y^{3}
 + ( - 3\,x + 2)\,y^{2} + ( - 2\,x^{2} + 3\,x)\,y + x^{3})^{2}
\]
\item $\Si_{red}=[C_{3,7},A_5,A_2,2A_1]$ with $C=B_1+B_1'+B_4$:
\[
\mathrm{f}(x, \,y) := ( - y^{2} + y - x^{2})^{3} + ( - 
{\displaystyle \frac {9}{2}} \,y^{3} + (3\,x + {\displaystyle 
\frac {9}{2}} )\,y^{2} - 3\,x\,y - x^{3})^{2}
\]
\item $\Si_{red}=[3A_2,C_{3,8},A_1]_1$, $C=B_1+B_5$, $\Si(B_5)=[3A_2,A_3,A_1]$:
\[
 \mathrm{f}(x, \,y) := (y^{2} + ( - 2\,x + 1)\,y - x^{2})^{3} + 
(y^{3} + ( - {\displaystyle \frac {3}{2}} \,x + 3)\,y^{2} + ( - 3
\,x^{2} - {\displaystyle \frac {3}{2}} \,x)\,y + x^{3})^{2}
\]
\item $\Si_{red}=[3A_2,C_{3,8},A_1]_2$, $C=B_1+B_5$, $\Si(B_5)=[3A_2,A_5]$:
\[
\mathrm{fff}(x, \,y) := (y^{2} + ( - x + 1)\,y - x^{2})^{3} + (y
^{3} + (x + 1)\,y^{2} + 3\,x\,y + x^{3})^{2}
\]
\item  $C=B_1+B_2+B_3$, $B_2\cap B_3=A_3+A_5+A_1$ and 
$\Si_{red}=[C_{3,8},A_5,A_2,A_1]_1$:

$$
\mathrm{f}(x, \,y) := ( - y^{2} + (x + {\displaystyle \frac {3}{
4}} )\,y - x^{2})^{3} + (y^{3} + ( - {\displaystyle \frac {3}{2}
} \,x + {\displaystyle \frac {3}{4}} )\,y^{2} - {\displaystyle 
\frac {9}{8}} \,y\,x + x^{3})^{2}
$$

\item  $C=B_1+B_2+B_3$, $B_2\cap B_3=A_5+A_5$ and 
$\Si_{red}=[C_{3,8},A_5,A_2,A_1]_2$: 
\[
 \mathrm{f}(x, \,y) := ( - {\displaystyle \frac {1}{16}} \,y^{2}
 + (x - 3)\,y - x^{2})^{3} + ({\displaystyle \frac {1}{64}} \,y^{
3} - {\displaystyle \frac {3}{8}} \,x\,y^{2} + (3\,x^{2} - 9\,x)
\,y + x^{3})^{2}\]
\item $C=B_1+B_5$, $\Si(B_5)=[A_8,A_3]$ and $\Si_{red}=[A_8,C_{3,8}]$:
\[
 \mathrm{f}(x,y) := ( - y^{2} + ( - 3\,x + 1)\,y - x^{2})^{3} + ( - 
{\displaystyle \frac {27}{8}} \,y^{3} + ( - {\displaystyle 
\frac {69}{8}} \,x + {\displaystyle \frac {27}{8}} )\,y^{2} + (
{\displaystyle \frac {9}{8}} \,x^{2} - {\displaystyle \frac {3}{2
}} \,x)\,y + x^{3})^{2}
\]
\item $\Si_{red}=[C_{3,9},E_6,A_1]$ with $C=B_1+B_5$:
\[
\mathrm{f}(x, \,y) := ( - y^{2} + y - x^{2})^{3} + (y^{2}\,x + 
( - 2\,x^{2} - x)\,y + x^{3})^{2}
\]
\item $[C_{3,12},A_2,A_1]_1 \mapright{u\to 1} [C_{3,12},A_2,A_1]_3$:
\[
\mathrm{f}(x, \,y) := ( - u\,y^{2} + y - x^{2})^{3} + (y^{3} + 
y^{2} + ( - x^{2} + x)\,y - x^{3})^{2}
\]
and $[C_{3,12},A_2,A_1]_2 \mapright{u\to 1} [C_{3,12},A_2,A_1]_3$:
\[
\mathrm{f}(x, \,y) := ( - y^{2} + u\,y - u\,x^{2})^{3} + (y^{3}
 + y^{2} + ( - x^{2} + x)\,y - x^{3})^{2}
\]
\item $\Si_{red}=[C_{6,6},2A_2,A_1]_1$ with $C=B_1+B_5$ $\mapright{s=0}$
$[C_{6,6},2A_2,2A_1]_1$ with $C=B_1+B_1'+B_4$:
\[
\mathrm{f}(x, \,y) := (y^{2} + x\,y - x^{2})^{3} + (y^{3} + ((
 - 7 + s)\,x + 2)\,y^{2} + (x^{2} - x)\,y + x^{3})^{2}
\]
\item $\Si_{red}=[C_{6,6},A_5,2A_1]$ with $C=B_1+B_1'+B_2+B_2'$:
\[
\mathrm{f}(x, \,y) :=  - y^{6} + (y^{3} + x\,y^{2} + (x^{2} + x
)\,y + 2\,x^{3} + x^{2})^{2}
\]

\item $C=B_3+B_3'$ with $[C_{6,12}]$  for $s\ne 1$ and 
$C=B_3+B_2+B_1$ with $[C_{6,12},A_1]$ for $s=1$ are given by:
\[
 \mathrm{f}(x, \,y) :=  - y^{6} + (y^{3} + (x + 1)\,y^{2} + (x^{
2} + x)\,y + 2\,s\,x^{3})^{2}
\]
\item $\Si_{red}=[C_{3,15},A_1]$ with $C=B_1+B_5$:
\[
 \mathrm{f}(x,y)=((x+1)y-x^2)^3+(y^3+(x+1)y^2+x\,y-x^3)^2
\]
\item $\Si_{red}=[B_{3,12}]$:
\begin{multline*}
\mathrm{f}(x, \,y) := ( - y^{2} + y - x^{2})^{3} + (y^{3} - 3\,
y^{2} + 3\,y\,x^{2})^{2}
\\
=( - x^{2} + 3\,y^{2} + 2\,y^{2}\,\sqrt{3} + y)\,(x^{2} - y - 3\,y
^{2} + 2\,y^{2}\,\sqrt{3})\,(x^{2} - y)
\end{multline*}

\item $C=B_1+B_5,\, \Si_{red}=[D_{4,7},2A_2]\mapright{s=0}
C=B_1+B_2+B_3,\,[D_{4,7},A_5]$:
\[
\mathrm{f}(x, \,y) := (s\,x\,y - x^{2})^{3} + ( - y^{3} + (x + 1
)\,y^{2} + y\,x^{2} + x^{3})^{2}
\]

\end{enumerate}

 \subsection{Examples II, Non-simple singularities}
 We give explicit examples for some configurations.
 The normal form of $B_{3,6}$-singulaity at $O$ with 
 $y=0$ a linear component is given as 
\begin{eqnarray*}
&\mathit{f_2} := \mathit{a_{02}}\,y^{2} + (\mathit{a_{11}}\,x + \mathit{
a_{01}})\,y - t^{2}\,x^{2}
\\
&\mathit{f_3} := \mathit{b_{03}}\,y^{3} + (\mathit{b12}\,x + \mathit{
b_{02}})\,y^{2} + (\mathit{b_{21}}\,x^{2} + \mathit{b_{11}}\,x)\,y + t^{3}
\,x^{3}
\end{eqnarray*}
The normal form of torus decomposition with $C=B_2+B_4$ and $B_{3,6}$ at
$O$
where $B_2$ is defined by $y-x^2=0$ is given
\begin{multline*}
\mathrm{f}(x, \,y) := ( - t^{2}\,y^{2} + (\mathit{a_{11}}\,x + 
\mathit{a_{01}})\,y - {\displaystyle \frac {1}{4}} \,{\displaystyle 
\frac {(4\,t^{2}\,\mathit{a_{01}} + \mathit{a_{11}}^{2})\,x^{2}}{t^{2}}
} )^{3} + {\displaystyle \frac {1}{64}} ( - 8\,t^{6}\,y^{3} \\
\mbox{} + 12\,y^{2}\,t^{4}\,\mathit{a_{11}}\,x - 8\,y^{2}\,t^{3}\,
\mathit{b_{02}} - 6\,x^{2}\,y\,t^{2}\,\mathit{a_{11}}^{2} + 8\,x^{2}\,y
\,t^{3}\,\mathit{b_{02}} - 8\,x\,y\,t^{3}\,\mathit{b_{11}} \\
\mbox{} + x^{3}\,\mathit{a_{11}}^{3} + 8\,x^{3}\,t^{3}\,\mathit{b_{11}}
)^{2}/t^{6} 
\end{multline*}
In the following examples, those with  a line or conic component can be easily
derived from
the above normal forms.
 \begin{enumerate}
\item $B_{3,6}\in \Si_{in}$:
$\Si_{red}=[4A_2,B_{36},A_1]$, $C=B_1+B_5,\, \Si(B_5)=[4A_2,A_3]$:

\begin{multline*}
\mathrm{f}(x,y)=
 (y^{2} + ( - 5\,x + 1)\,y - x^{2})^{3} + (
 - {\displaystyle \frac {3}{4}} \,y^{3} + ({\displaystyle \frac {
15}{4}} \,x + {\displaystyle \frac {3}{4}} )\,y^{2} + ( - 15\,x^{
2} + {\displaystyle \frac {9}{4}} \,x)\,y + x^{3})^{2}
\end{multline*}
\item  $\Si_{red}=[A_2,A_5,B_{3,6},2A_1]$, $C=B_1+B_1'+B_4$:
\[
 \mathrm{f}(x,y)
= ( - y^{2} + y - 4\,x^{2})^{3} + (y^{3} + ( - 4\,x - 1)\,y^{2} + 4\,y\,x - 8\,x^{3})^{2}
\]
\item  $\Si_{red}=[A_2,A_5,B_{3,6},2A_1]$, $C=B_3+B_2+B_1$:
\[
 \mathrm{f}(x,y)
= ( - y^{2} + ( - 3\,x + 1)\,y - x^{2})^{3}
 + ({\displaystyle \frac {2}{3}} \,y^{3} + ({\displaystyle 
\frac {13}{3}} \,x - {\displaystyle \frac {2}{3}} )\,y^{2} + (5\,
x^{2} - {\displaystyle \frac {7}{3}} \,x)\,y + x^{3})^{2}
\]
\item  $[A_8,B_{3,6},A_1]$, $C=B_1+B_5$ and $\Si(B_5)=[A_8,A_3]$:
\[
 \mathrm{f}(x,y)=
( - y^{2} + ( - 3\,x + 1)\,y - x^{2})^{3}
 + ( - y^{3} + ({\displaystyle \frac {2}{3}} \,x + 1)\,y^{2} + (
10\,x^{2} - {\displaystyle \frac {11}{3}} \,x)\,y + x^{3})^{2}
\]
\item $\Si_{red}=[C_{3,7},A_8,A_1]$ with $C=B_1+B_5$:
\[
\mathrm{f}(x, \,y) := ( - y^{2} + y - x^{2})^{3} + ( - 2\,y^{3}
 + ( - 3\,x + 2)\,y^{2} + ( - 2\,x^{2} + 3\,x)\,y + x^{3})^{2}
\]
\item $\Si_{red}=[C_{3,7},A_5,A_2,2A_1]$ with $C=B_1+B_1'+B_4$:
\[
\mathrm{f}(x, \,y) := ( - y^{2} + y - x^{2})^{3} + ( - 
{\displaystyle \frac {9}{2}} \,y^{3} + (3\,x + {\displaystyle 
\frac {9}{2}} )\,y^{2} - 3\,x\,y - x^{3})^{2}
\]
\item $\Si_{red}=[3A_2,C_{3,8},A_1]_1$, $C=B_1+B_5$, $\Si(B_5)=[3A_2,A_3,A_1]$:
\[
 \mathrm{f}(x, \,y) := (y^{2} + ( - 2\,x + 1)\,y - x^{2})^{3} + 
(y^{3} + ( - {\displaystyle \frac {3}{2}} \,x + 3)\,y^{2} + ( - 3
\,x^{2} - {\displaystyle \frac {3}{2}} \,x)\,y + x^{3})^{2}
\]
\item $\Si_{red}=[3A_2,C_{3,8},A_1]_2$, $C=B_1+B_5$, $\Si(B_5)=[3A_2,A_5]$:
\[
\mathrm{fff}(x, \,y) := (y^{2} + ( - x + 1)\,y - x^{2})^{3} + (y
^{3} + (x + 1)\,y^{2} + 3\,x\,y + x^{3})^{2}
\]
\item  $C=B_1+B_2+B_3$, $B_2\cap B_3=A_3+A_5+A_1$ and 
$\Si_{red}=[C_{3,8},A_5,A_2,A_1]_1$:

$$
\mathrm{f}(x, \,y) := ( - y^{2} + (x + {\displaystyle \frac {3}{
4}} )\,y - x^{2})^{3} + (y^{3} + ( - {\displaystyle \frac {3}{2}
} \,x + {\displaystyle \frac {3}{4}} )\,y^{2} - {\displaystyle 
\frac {9}{8}} \,y\,x + x^{3})^{2}
$$

\item  $C=B_1+B_2+B_3$, $B_2\cap B_3=A_5+A_5$ and 
$\Si_{red}=[C_{3,8},A_5,A_2,A_1]_2$: 
\[
 \mathrm{f}(x, \,y) := ( - {\displaystyle \frac {1}{16}} \,y^{2}
 + (x - 3)\,y - x^{2})^{3} + ({\displaystyle \frac {1}{64}} \,y^{
3} - {\displaystyle \frac {3}{8}} \,x\,y^{2} + (3\,x^{2} - 9\,x)
\,y + x^{3})^{2}\]
\item $C=B_1+B_5$, $\Si(B_5)=[A_8,A_3]$ and $\Si_{red}=[A_8,C_{3,8}]$:
\[
 \mathrm{f}(x,y) := ( - y^{2} + ( - 3\,x + 1)\,y - x^{2})^{3} + ( - 
{\displaystyle \frac {27}{8}} \,y^{3} + ( - {\displaystyle 
\frac {69}{8}} \,x + {\displaystyle \frac {27}{8}} )\,y^{2} + (
{\displaystyle \frac {9}{8}} \,x^{2} - {\displaystyle \frac {3}{2
}} \,x)\,y + x^{3})^{2}
\]
\item $\Si_{red}=[C_{3,9},E_6,A_1]$ with $C=B_1+B_5$:
\[
\mathrm{f}(x, \,y) := ( - y^{2} + y - x^{2})^{3} + (y^{2}\,x + 
( - 2\,x^{2} - x)\,y + x^{3})^{2}
\]
\item $[C_{3,12},A_2,A_1]_1 \mapright{u\to 1} [C_{3,12},A_2,A_1]_3$:
\[
\mathrm{f}(x, \,y) := ( - u\,y^{2} + y - x^{2})^{3} + (y^{3} + 
y^{2} + ( - x^{2} + x)\,y - x^{3})^{2}
\]
and $[C_{3,12},A_2,A_1]_2 \mapright{u\to 1} [C_{3,12},A_2,A_1]_3$:
\[
\mathrm{f}(x, \,y) := ( - y^{2} + u\,y - u\,x^{2})^{3} + (y^{3}
 + y^{2} + ( - x^{2} + x)\,y - x^{3})^{2}
\]
\item $\Si_{red}=[C_{6,6},2A_2,A_1]_1$ with $C=B_1+B_5$ $\mapright{s=0}$
$[C_{6,6},2A_2,2A_1]_1$ with $C=B_1+B_1'+B_4$:
\[
\mathrm{f}(x, \,y) := (y^{2} + x\,y - x^{2})^{3} + (y^{3} + ((
 - 7 + s)\,x + 2)\,y^{2} + (x^{2} - x)\,y + x^{3})^{2}
\]
\item $\Si_{red}=[C_{6,6},A_5,2A_1]$ with $C=B_1+B_1'+B_2+B_2'$:
\[
\mathrm{f}(x, \,y) :=  - y^{6} + (y^{3} + x\,y^{2} + (x^{2} + x
)\,y + 2\,x^{3} + x^{2})^{2}
\]

\item $C=B_3+B_3'$ with $[C_{6,12}]$  for $s\ne 1$ and 
$C=B_3+B_2+B_1$ with $[C_{6,12},A_1]$ for $s=1$ are given by:
\[
 \mathrm{f}(x, \,y) :=  - y^{6} + (y^{3} + (x + 1)\,y^{2} + (x^{
2} + x)\,y + 2\,s\,x^{3})^{2}
\]
\item $\Si_{red}=[C_{3,15},A_1]$ with $C=B_1+B_5$:
\[
 \mathrm{f}(x,y)=((x+1)y-x^2)^3+(y^3+(x+1)y^2+x\,y-x^3)^2
\]
\item $\Si_{red}=[B_{3,12}]$:
\begin{multline*}
\mathrm{f}(x, \,y) := ( - y^{2} + y - x^{2})^{3} + (y^{3} - 3\,
y^{2} + 3\,y\,x^{2})^{2}
\\
=( - x^{2} + 3\,y^{2} + 2\,y^{2}\,\sqrt{3} + y)\,(x^{2} - y - 3\,y
^{2} + 2\,y^{2}\,\sqrt{3})\,(x^{2} - y)
\end{multline*}

\item $C=B_1+B_5,\, \Si_{red}=[D_{4,7},2A_2]\mapright{s=0}
C=B_1+B_2+B_3,\,[D_{4,7},A_5]$:
\[
\mathrm{f}(x, \,y) := (s\,x\,y - x^{2})^{3} + ( - y^{3} + (x + 1
)\,y^{2} + y\,x^{2} + x^{3})^{2}
\]

\end{enumerate}

\subsection{Non-trivial degenerations}
We give some non-trivial degenerations. 
 We do not give the whole moduli description but it can
   be computed as in \cite{Oka-Pho2}.
   \begin{enumerate}
\item  (i) $\Si_{red}=[A_5,4A_2,A_3,2A_1]_2$, $C=B_1+B_5$ and
       $\Si(B_5)=[4A_2,A_3]$ can not degenerate into any further simple
       configuration, but it 
degenerates into $C_0=B_1+B_1'+B_4$ with
       $\Si_{red}=[A_5,C_{3,7},A2,2A_1]$ as $t\to 0$.
For $t=0$, $A_5$ and $C_{3,7}$ are at $O$ and $(0,1)$.
\begin{multline*}
\mathrm{f}_t(x, \,y) 
=(( - 1 - t^{2})\,y^{2} + ( - 4\,t\,x + 1)\,y - x^{2})^{3} \\
\mbox{} + (( - t^{3} - {\displaystyle \frac {3}{2}} \,t + 1)\,y^{
3} + (( - 3 - 6\,t^{2})\,x + {\displaystyle \frac {3}{2}} \,t - 2
)\,y^{2} + ( - {\displaystyle \frac {15}{2}} \,t\,x^{2} + 3\,x + 
1)\,y + x^{3})^{2} 
\end{multline*}
\item $\Si_{red}=[A_5,4A_2,A_3,2A_1]_3$, $C=B_1+B_5$,
       $\Si(B_5)=[4A_4,2A_1]$
degenerates into $C_0=B_1+B_5$, with
       $\Si_{red}=[A_5,E_6,A_3,2A_2,A_1]^{mr}$ 
when        $t\to 0$.
\begin{multline*}
\mathrm{f}_t(x, \,y) =( - y^{2} - y^{2}\,t^{2} - 3\,y\,x - 4\,y\,x\,t + y - x^{2})^{3}
 + {\displaystyle \frac {1}{64}} (8\,y^{3}\,t^{4} + 48\,y^{2}\,x
\,t^{3} + 8\,y^{3}\,t^{3} + 60\,y\,x^{2}\,t^{2} \\
\mbox{} - 12\,y^{2}\,t^{2} + 84\,y^{2}\,x\,t^{2} + 12\,y^{3}\,t^{
2} - 12\,y^{2}\,t + 12\,y^{3}\,t + 132\,y\,x^{2}\,t + 60\,y^{2}\,
x\,t - 24\,y\,x\,t \\
\mbox{} - 8\,x^{3}\,t - 6\,y^{2} - 8\,x^{3} - 24\,y\,x + 72\,y\,x
^{2} + 24\,y^{2}\,x + 3\,y + 3\,y^{3})^{2} \left/ {\vrule 
height0.44em width0em depth0.44em} \right. \!  \! (1 + t)^{2} 
\end{multline*}
\item $\Si_{red}=[A_5,4A_2,D_5]_2$, $C=B_1+B_5$ 
degenerates into
$\Si_{red}=[2A_5,2A_2,D_5]_1^{mr}$:
For $t=0$, two $A_5$'s are at $\{(0,0), (0,1)\}$.
\begin{multline*}
\mathrm{f}_2(x,y)=
 - y^{2} + (( - {\displaystyle \frac {25}{16}}  + {\displaystyle 
\frac {1}{16}} \,t)\,x + 1)\,y - x^{2}\\
\mathrm{f}_3(x,y)= ( - y^{2} + (( - {\displaystyle \frac {25
}{16}}  + {\displaystyle \frac {1}{16}} \,t)\,x + 1)\,y - x^{2})
^{3} + {\displaystyle \frac {1}{4294967296}} (6428160\,y^{3} + 
65536\,x^{3}\,t^{2} \\
\mbox{} + 65536\,y\,t^{2} - 1179648\,y\,t - 11736576\,y^{2} + 
1552896\,y^{2}\,t - 41472\,y^{2}\,t^{2} \\
\mbox{} + 1536\,y^{2}\,t^{3} - 36\,y^{3}\,t^{4} - 2586\,y^{3}\,t
^{3} + y^{3}\,t^{5} - 491103\,y^{3}\,t + 66204\,y^{3}\,t^{2} + 
5308416\,y \\
\mbox{} + 4128768\,y\,x\,t - 229376\,y\,x\,t^{2} - 6144\,y\,x^{2}
\,t^{3} + 329728\,y\,x^{2}\,t^{2} - 4442112\,y\,x^{2}\,t \\
\mbox{} + 96\,y^{2}\,x\,t^{4} - 4992\,y^{2}\,x\,t^{3} + 186432\,y
^{2}\,x\,t^{2} - 3602304\,y^{2}\,x\,t - 18579456\,y\,x \\
\mbox{} + 5308416\,x^{3} + 20329056\,y^{2}\,x - 1179648\,x^{3}\,t
 + 17750016\,x^{2}\,y)^{2} \left/ {\vrule 
height0.44em width0em depth0.44em} \right. \!  \! ( - 9 + t)^{4}
\end{multline*}
\item   $[2A_5,2A_2,4A_1], \,C=B_4+B_1+B_1'\to [3A_5,4A_1]^{mr},\, B_3+B_1+B_1'+B_1''$:
\begin{multline*}
\mathrm{f}(x,y,a)=(y\,x\,a^{2} - x^{2}\,a + a\,x - 4\,y\,x\,
a - y^{2})^{3} + ({\displaystyle \frac {1}{2}} \,y^{2}\,x\,a^{3}
 - {\displaystyle \frac {1}{2}} \,a^{3}\,x^{2}\,y - 
{\displaystyle \frac {1}{2}} \,x^{2}\,a^{2} + {\displaystyle 
\frac {3}{2}} \,y\,x^{2}\,a^{2} \\
\mbox{} - {\displaystyle \frac {7}{2}} \,y^{2}\,x\,a^{2} + 
{\displaystyle \frac {1}{2}} \,x^{3}\,a^{2} + y\,x\,a^{2} - 
{\displaystyle \frac {1}{2}} \,x^{2}\,a + 2\,y\,x^{2}\,a - 3\,y\,
x\,a + {\displaystyle \frac {1}{2}} \,a\,x + 7\,y^{2}\,x\,a - 4\,
y^{2}\,x \\
\mbox{} + {\displaystyle \frac {1}{2}} \,x - 3\,y\,x^{2} - y\,x
 - y^{3} - {\displaystyle \frac {1}{2}} \,x^{3})^{2} 
\end{multline*}
The curve $f(x,y,0)$
has three line components and a nodal cubic
component.
\[
f(x,y,0)={\displaystyle \frac {1}{4}} \,x\,(x + 1 + 2\,y)\,(x - 1 + 4\,y)
\,(4\,y^{3} + 8\,y^{2}\,x + 6\,y\,x^{2} + 2\,y\,x + x^{3} - x)
\]
\item   $C_u=B_2+B_4$, $\Si_{red}=[2A_5,3A_2,2A_1]_2\to$
 $ [3A_5,A_2,2A_1]^{mr}$ with $C_0=B_1+B_2+B_3$.
\begin{multline*}
\mathrm{f}_2(x,y,u)= - y^{2} + y^{2}\,\mathit{t_1}\,u - {\displaystyle \frac {1}{2}} 
\,y^{2}\,u^{2} - x\,y\,\mathit{t_1} + x\,y\,u - {\displaystyle 
\frac {1}{4}} \,x^{2}\,\mathit{t_1}^{2} + {\displaystyle \frac {1
}{2}} \,x^{2}\,\mathit{t_1}\,u - {\displaystyle \frac {1}{4}} \,x
^{2}\,u^{2} - x\,\mathit{t_1}\,u + {\displaystyle \frac {1}{2}} \,
x\,u^{2}\\
\mathrm{f}_3(x,y,u)=y^{3} + {\displaystyle \frac {3}{8}} \,x^{3}\,\mathit{t_1}\,u^{2}
 + {\displaystyle \frac {3}{2}} \,y^{3}\,u^{2} + {\displaystyle 
\frac {9}{4}} \,x^{2}\,\mathit{t_1}^{2}\,u - {\displaystyle 
\frac {1}{2}} \,y^{2}\,x\,u^{3} + {\displaystyle \frac {3}{4}} \,
x^{2}\,y\,u^{2} - {\displaystyle \frac {3}{2}} \,x\,y\,u^{2} - 
{\displaystyle \frac {9}{2}} \,y^{3}\,\mathit{t_1}\,u \\- 
{\displaystyle \frac {3}{2}} \,x^{2}\,y\,\mathit{t_1}\,u
 + {\displaystyle \frac {9}{2}} \,x\,y\,\mathit{t_1}\,u + 
{\displaystyle \frac {15}{8}} \,y^{2}\,x\,\mathit{t_1}\,u^{2} - 
{\displaystyle \frac {9}{4}} \,y^{2}\,x\,\mathit{t_1}^{2}\,u - 
{\displaystyle \frac {9}{4}} \,u^{2}\,y^{2}\,\mathit{t_1} + 3\,u\,
y^{2}\,\mathit{t_1}^{2} - 3\,u\,x\,\mathit{t_1}^{2} \\
\mbox{} - {\displaystyle \frac {1}{2}} \,x\,u^{3} - 
{\displaystyle \frac {1}{8}} \,x^{3}\,u^{3} + {\displaystyle 
\frac {1}{2}} \,x^{2}\,u^{3} + {\displaystyle \frac {1}{2}} \,y^{
2}\,u^{3} - {\displaystyle \frac {3}{8}} \,x^{3}\,\mathit{t_1}^{2}
\,u + {\displaystyle \frac {9}{4}} \,u^{2}\,x\,\mathit{t_1} - 
{\displaystyle \frac {15}{8}} \,x^{2}\,\mathit{t_1}\,u^{2} - y^{2}
\,\mathit{t_1}^{3} \\
\mbox{} + y^{2}\,x\,\mathit{t_1}^{3} + x\,\mathit{t_1}^{3} + 
{\displaystyle \frac {1}{8}} \,x^{3}\,\mathit{t_1}^{3} + 3\,y^{3}
\,\mathit{t_1}^{2} - x^{2}\,\mathit{t_1}^{3} + {\displaystyle 
\frac {3}{2}} \,y^{2}\,x\,\mathit{t_1} - 3\,x\,y\,\mathit{t_1}^{2}
 + {\displaystyle \frac {3}{4}} \,x^{2}\,y\,\mathit{t_1}^{2} \\
\mbox{} - {\displaystyle \frac {3}{2}} \,y^{2}\,x\,u 
\end{multline*}
with $\mathit{t_1} := u + {\displaystyle \frac {1}{2}} \,\sqrt{u^{2} - 6}$.
\begin{multline*}
f(x,y,0)={\displaystyle \frac {3}{32}} (9\,x^{3} - 36\,x^{2} + 36\,x - 18
\,x^{2}\,y\,\sqrt{-6} + 36\,x\,y\,\sqrt{-6} - 36\,y^{2}\,x - 36\,
y^{2} - 20\,y^{3}\,\sqrt{-6}) \\
(x - 1 - \sqrt{-6}\,y)\,(x - y^{2}) 
\end{multline*}
\item $[3A_5,A_3]_1, \,C=B_1+B_5\mapright{u\to 1} [3A_5,D_4],\, 
C=B_1+B_1+B_4$:
\begin{multline*}
\mathrm{f}(x, \,y) := ( - y^{2} - 3\,x\,y - x^{2} + x)^{3} + 
{\displaystyle \frac {1}{16}} ( - 4\,y^{3}\,u - 14\,x\,y^{2}\,u
 + x\,y^{2} + 4\,x\,y^{2}\,u^{2} - 10\,x^{2}\,y\,u \\
\mbox{} + 3\,x^{2}\,y + 4\,x^{2}\,y\,u^{2} + 4\,x\,y\,u - 4\,x\,y
\,u^{2} - 2\,u\,x^{3} + x^{3} + u^{2}\,x^{3} + 2\,x^{2}\,u - x^{2
} - 2\,x^{2}\,u^{2} \\
\mbox{} + x\,u^{2})^{2}/u^{2} 
\end{multline*}
\item $[3A_5,A_3]_2,\, C=B_2+B_4\mapright{} [3A_5,D_4]_1^{mr}\, (a=1/4),\, [3A_5,D_4]_2^{mr}\,(a=-1/12)$:
\begin{multline*}
\mathrm{f}(x, \,y) := ( - y^{2}\,a - {\displaystyle \frac {1}{2}
} \,y^{2} + {\displaystyle \frac {1}{2}} \,y - a\,x^{2} + a)^{3}
 \\
\mbox{} + ({\displaystyle \frac {5}{4}} \,y^{3}\,a + 
{\displaystyle \frac {3}{8}} \,y^{3} - {\displaystyle \frac {1}{4
}} \,y^{2}\,a - {\displaystyle \frac {1}{2}} \,y^{2} + 
{\displaystyle \frac {5}{4}} \,y\,a\,x^{2} - {\displaystyle 
\frac {5}{4}} \,y\,a + {\displaystyle \frac {1}{8}} \,y - 
{\displaystyle \frac {1}{4}} \,a\,x^{2} + {\displaystyle \frac {1
}{4}} \,a)^{2} 
\end{multline*}
 \item A degeneration $[A_{11},2A_2,A_3] \, s\ne 1\to [A_{11},2A_2,D_4],\, s=1$ with 
$C=B_2+B_4$. 
\begin{multline*}
\mathrm{f}(x, \,y) := ( - y^{2} + y - {\displaystyle \frac {8}{3
}} \,{\displaystyle \frac {x^{2}}{s}} )^{3} + (y^{3} - 
{\displaystyle \frac {5}{4}} \,y^{2} + ({\displaystyle \frac {10
}{3}} \,{\displaystyle \frac {x^{2}}{s}}  + x + {\displaystyle 
\frac {3}{8}} \,s)\,y - {\displaystyle \frac {8}{3}} \,
{\displaystyle \frac {x^{3}}{s}}  - x^{2})^{2}
\end{multline*}
\end{enumerate}
\section{Appendix}
{\bf Proof of Proposition \ref{linear-torus}}.
The proof is computational. 
\nl
1.  Assume that $C$ is a sextics of linear type with $3A_5$.
This is a special case of the result of Tokunaga, \cite{Tokunaga-torus}.
We assume that $L$ is defined by $y=0$.
We start from the generic polynomial of
degree 6:
\[
\mathrm{f}(x,y)=\sum_{i+j\le 6}a_{i,j} x^i y^j
\]
By the action of $PGL(3,\bfC)$, we may assume that three $A_5$'s are at
$\{P_1=(-1,0),P_2=(0,0),P_3= (1,0)\}$
and the tangent cones  at $P_1, P_3$ are given by 
$x=\pm 1$. These condition says

(1) $f(x,0)=x^2(x^2-1)^2$ and 

(2) $ f_x(P_i)=f_y(P_i)=0$ for $i=1,2,3$ and

(3) $f_{x,y}(P_i)=f_{y,y}(P_i)=0, \, i=1,3$.

Eliminating coefficients from $f(x,y)$ using these equations, then we 
elliminate further coefficients using the assumptions
$(C,P_i)\cong A_5$. Then we get a normal form of sextics with 
four parameters $a_{04},a_{06},t,s$. Then we apply the  degeneration
method  to the family $fu:=f-u y^6$ (see \cite{Oka-Pho2}). 
Finally we find that $f(x,y)$ has a torus expression:
\begin{multline*}
\mathrm{f}(x,y)=\tau\,y^{6} + 
( - t\,y^{2} - s\,y^{2} - y\,x^{2} + x\,t\,y^{2}\\
 + y - x\,s\,y^{2
} + x^{3} - x - y^{3}\,s\,t + {\displaystyle \frac {1}{2}} \,y^{3
}\,\mathit{a_{04}} - {\displaystyle \frac {1}{2}} \,y^{3}\,s^{2} - 
{\displaystyle \frac {1}{2}} \,y^{3}\,t^{2})^{2} \\
\tau=(\mathit{a_{06}} + {\displaystyle \frac {1}{2}} \,s^{2}\,\mathit{a_{04}
} - t^{3}\,s - s^{3}\,t - {\displaystyle \frac {3}{2}} \,s^{2}\,t
^{2} - {\displaystyle \frac {1}{4}} \,t^{4} + t\,\mathit{a_{04}}\,s
 - {\displaystyle \frac {1}{4}} \,s^{4} - {\displaystyle \frac {1
}{4}} \,\mathit{a_{04}}^{2} + {\displaystyle \frac {1}{2}} \,
\mathit{a_{04}}\,t^{2})
\end{multline*}
The proof for the cases $\Si(C)\supset A_{11}+A_5, A_{17}$ are similar.
In the case of $A_{11}+A_5$, we assume that 
$A_{11}$ is at $O$ and  $A_5$ is at $(1,0)$.
We assume that the tangent cone at $(1,0)$ is $x=1$.
In the case $A_{17}$, we assume that $A_{17}$ is at $O$.
For the condition of $(C,O)\cong A_{11}$ (respectively for 
$\cong A_{17}$) at the origin,
we assume that 
the normal form is given by the change of coordinates
$y_1:= y+\sum_{i=2}^5 t_i x^i$ (resp. $y_1:= y+\sum_{i=2}^8 t_i
x^i$). Let 
$f'(x,y)=f(x,y_1)$ and let 
$c_0=f'(x,0)$ and $ c_1$ be the coefficients of $y$ in $f'(x,y)$.
As we assume that $f'(x,y)=a y^2+b y x^6+ c x^{12} +\text{(higher terms)s}$
(resp. $f'(x,y)=a y^2+b y x^9+ c x^{18} +\text{(higher terms)s}$),
normal forms are
obtained by solving the equalities:
$ Coeff(c_0,x,j)=0$ for $j\le 11$ and $Coeff(c_1,x,k)=0$ for $k\le 5$
(resp. $ Coeff(c_0,x,j)=0$ for $j\le 17$ and $Coeff(c_1,x,k)=0$ for
$k\le 8$).
Then we find the torus expressions by the degeneration method. 
For the case $A_{11}+A_5$,
$f(x,y)=u_{11} y^6+f_3(x,y)^2$
 where
\begin{multline*}
\mathit{u_{11}} =  - {\displaystyle \frac {1}{4}} (4\,\mathit{a_{04}}
\,\mathit{t_2}^{9}\,\mathit{t_4}\,\mathit{t_3}^{2} + \mathit{t_2}^{4}
\,\mathit{t_4}^{4} - 2\,\mathit{a_{04}}\,\mathit{t_2}^{10}\,\mathit{t_4
}^{2} - 2\,\mathit{t_2}^{8}\,\mathit{a_{04}}\,\mathit{t_3}^{4} + 6\,
\mathit{t_2}^{2}\,\mathit{t_4}^{2}\,\mathit{t_3}^{4} \\
\mbox{} - 4\,\mathit{t_2}^{3}\,\mathit{t_4}^{3}\,\mathit{t_3}^{2} - 
4\,\mathit{t_3}^{6}\,\mathit{t_4}\,\mathit{t_2} + \mathit{a_{04}}^{2}\,
\mathit{t_2}^{16} + \mathit{t_3}^{8} - 4\,\mathit{a_{06}}\,\mathit{t_2}
^{14})/\mathit{t_2}^{14}\\
\mathrm{f_3}(x, \,y) = {\displaystyle \frac {1}{2}} (6\,
\mathit{t_2}^{3}\,y^{2}\,\mathit{t_3}\,x\,\mathit{t_4} - 2\,\mathit{
t_2}^{7}\,x^{2} - 2\,\mathit{t_2}^{4}\,y^{2}\,\mathit{t_4} - 2\,
\mathit{t_2}^{6}\,y\,x + 2\,\mathit{t_2}^{3}\,y^{2}\,\mathit{t_3}^{2
} \\
\mbox{} - 2\,\mathit{t_2}^{3}\,y^{2}\,\mathit{t_3}^{2}\,x + 2\,
\mathit{t_2}^{4}\,y^{2}\,x\,\mathit{t_4} - 4\,\mathit{t_2}^{2}\,y^{2
}\,\mathit{t_3}^{3}\,x + 2\,\mathit{t_2}^{6}\,y - 2\,\mathit{t_2}^{4
}\,y^{2}\,x\,\mathit{t_5} - y^{3}\,\mathit{t_2}^{2}\,\mathit{t_4}^{2
} \\
\mbox{} + 2\,\mathit{t_2}\,y^{3}\,\mathit{t_3}^{2}\,\mathit{t_4} - y
^{3}\,\mathit{t_3}^{4} + 2\,\mathit{t_2}^{7}\,x^{3} + \mathit{t_2}^{
8}\,y^{3}\,\mathit{a_{04}} - 2\,\mathit{t_2}^{5}\,y\,\mathit{t_3}\,x
 + 2\,\mathit{t_2}^{5}\,y\,\mathit{t_3}\,x^{2})/\mathit{t_2}^{7}\end{multline*}
and $f(x,y)=u_{17}y^6+h_3(x,y)^2$ for the case $A_{17}$ where
\begin{multline*}
\mathit{u_{17}}=  - {\displaystyle \frac {1}{4}} (24\,\mathit{t_6}
^{2}\,\mathit{t_3}^{6}\,\mathit{t_5}^{2}\,\mathit{t_4}^{2} - 24\,
\mathit{t_4}^{4}\,\mathit{t_6}^{2}\,\mathit{t_3}^{5}\,\mathit{t_5} + 
8\,\mathit{a_{04}}\,\mathit{t_3}^{13}\,\mathit{t_4}\,\mathit{t_6}\,
\mathit{t_5} - 8\,\mathit{t_4}\,\mathit{t_6}^{3}\,\mathit{t_3}^{7}\,
\mathit{t_5} \\
\mbox{} - 8\,\mathit{a_{04}}\,\mathit{t_3}^{12}\,\mathit{t_5}^{2}\,
\mathit{t_4}^{2} + 8\,\mathit{a_{04}}\,\mathit{t_4}^{4}\,\mathit{t_3}^{
11}\,\mathit{t_5} + \mathit{a_{04}}^{2}\,\mathit{t_3}^{20} - 2\,
\mathit{a_{04}}\,\mathit{t_4}^{6}\,\mathit{t_3}^{10} - 2\,\mathit{a_{04}}
\,\mathit{t_3}^{14}\,\mathit{t_6}^{2} \\
\mbox{} + 4\,\mathit{t_4}^{3}\,\mathit{t_6}^{3}\,\mathit{t_3}^{6} + 
6\,\mathit{t_4}^{6}\,\mathit{t_6}^{2}\,\mathit{t_3}^{4} - 4\,
\mathit{a_{04}}\,\mathit{t_3}^{12}\,\mathit{t_4}^{3}\,\mathit{t_6} + 
\mathit{t_3}^{8}\,\mathit{t_6}^{4} + 4\,\mathit{t_4}^{9}\,\mathit{t_6
}\,\mathit{t_3}^{2} - 4\,\mathit{a_{06}}\,\mathit{t_3}^{18} \\
\mbox{} + 24\,\mathit{t_4}^{8}\,\mathit{t_5}^{2}\,\mathit{t_3}^{2}
 - 8\,\mathit{t_4}^{10}\,\mathit{t_5}\,\mathit{t_3} + 48\,\mathit{t_4
}^{5}\,\mathit{t_6}\,\mathit{t_3}^{4}\,\mathit{t_5}^{2} - 24\,
\mathit{t_4}^{7}\,\mathit{t_6}\,\mathit{t_3}^{3}\,\mathit{t_5} + 
\mathit{t_4}^{12} \\
\mbox{} - 32\,\mathit{t_6}\,\mathit{t_3}^{5}\,\mathit{t_5}^{3}\,
\mathit{t_4}^{3} + 16\,\mathit{t_4}^{4}\,\mathit{t_5}^{4}\,\mathit{
t_3}^{4} - 32\,\mathit{t_4}^{6}\,\mathit{t_5}^{3}\,\mathit{t_3}^{3})/
\mathit{t_3}^{18}\\
\mathrm{h_3}(x, \,y) := {\displaystyle \frac {1}{2}} ( - 4\,
\mathit{t_5}^{2}\,\mathit{t_3}^{2}\,y^{3}\,\mathit{t_4}^{2} + 2\,
\mathit{t_5}^{2}\,\mathit{t_3}^{5}\,y^{2}\,x + 4\,\mathit{t_6}\,
\mathit{t_5}\,\mathit{t_3}^{3}\,y^{3}\,\mathit{t_4} - 10\,\mathit{t_5
}\,\mathit{t_3}^{4}\,y^{2}\,x\,\mathit{t_4}^{2} \\
\mbox{} + 4\,\mathit{t_5}\,\mathit{t_3}^{5}\,y^{2}\,\mathit{t_4} + 4
\,\mathit{t_5}\,\mathit{t_3}\,y^{3}\,\mathit{t_4}^{4} - 2\,\mathit{
t_5}\,\mathit{t_3}^{7}\,y\,x^{2} - 2\,\mathit{t_3}^{9}\,x^{3} + 2\,
\mathit{t_3}^{6}\,x^{2}\,\mathit{t_4}^{2}\,y \\
\mbox{} - 2\,\mathit{t_3}^{6}\,y^{2}\,x\,\mathit{t_7} + 6\,\mathit{
t_6}\,\mathit{t_3}^{5}\,y^{2}\,x\,\mathit{t_4} + 4\,\mathit{t_3}^{3}
\,x\,y^{2}\,\mathit{t_4}^{4} - 2\,\mathit{t_3}^{7}\,x\,y\,\mathit{
t_4} - y^{3}\,\mathit{t_6}^{2}\,\mathit{t_3}^{4} \\
\mbox{} - 2\,\mathit{t_6}\,\mathit{t_3}^{2}\,y^{3}\,\mathit{t_4}^{3}
 - y^{3}\,\mathit{t_4}^{6} - 2\,\mathit{t_6}\,\mathit{t_3}^{6}\,y^{2
} - 2\,\mathit{t_3}^{4}\,y^{2}\,\mathit{t_4}^{3} + 2\,\mathit{t_3}^{
8}\,y + y^{3}\,\mathit{t_3}^{10}\,\mathit{a_{04}})/\mathit{t_3}^{9} 
\end{multline*}

\noindent
2. Next we consider the case $C=B_3+B_3'$ and $C$ is a sextics of torus
type
and assume that $\Si_{in}(C)=[3A_5]$ or $[A_{11},A_5]$
or $[A_{17}]$ and let $P_1,P_2,P_3$ be the 
corresponding singular points.
(In the case of $[A_{11},A_5]$ or $[A_{17}]$, $P_2=P_3$ and
$P_1=P_2=P_3$ respectively.)
We show that  three $P_1,P_2,P_3$ must be colinear. We start from the
expression:
\begin{multline*}
\mathrm{f}(x, \,y) :=f_{31}(x,y) f_{32}(x,y)\\
= (\mathit{a_{03}}\,y^{3} + (\mathit{a_{12}}\,x + 
\mathit{a_{02}})\,y^{2} + (\mathit{a_{21}}\,x^{2} + \mathit{a_{11}}\,x + 
\mathit{a_{01}})\,y + \mathit{a_{30}}\,x^{3} 
+ \mathit{a_{20}}\,x^{2} + \mathit{a_{10}}\,x \mbox{} + \mathit{a_{00}})\\
(\mathit{b_{03}}\,y^{3} + (\mathit{b_{12}}\,x + 
\mathit{b_{02}})\,y^{2} + (\mathit{b_{21}}\,x^{2} + \mathit{b_{11}}\,x 
+ \mathit{b_{01}})\,y + \mathit{b_{30}}\,x^{3} + \mathit{b_{20}}\,x^{2} 
\mbox{} + \mathit{b_{10}}\,x + \mathit{b_{00}}) 
\end{multline*}
Assume that $C$ is defined by $f_2(x,y)^3+f_3(x,y)^2=0$
and let $C_2$ and $C_3$ be the conic and the cubic defined by $f_2=0$ 
and $f_3=0$. Let $P_i\in C$ be an inner singularity.
Recall that by \cite {Pho} we have the equivalence

($(\star)$ $(C,P_i)\cong A_{6j-1} \iff I(C_2,C_3;P_i)=2j$ and 
$C_3$ is smooth at $P_i$ for $j=1,2,3$.

In particular, if $C_2$ is smooth at $P_i$, $(\star)$ implies that $C_2$ and $C_3$ are 
tangent at $P_i$.
Again by an easy computation, we can see that there are no cases when
$P_1,P_2,P_3$ are not colinear. We give a recipe of the computation.
Assuming that $P_1,P_2,P_3$ are not  colinear.
To each $P_i$, we associate its  tangent cone direction $\ell_i$
and we identify $\ell_i$ and  a line in $\bfP^2$.

There are two cases.

(a) $C_2$ is a smooth conic, or 

(b) $C_2$ is a union of two distinct
lines
$L_1, L_2$.

In the case of (a), we may assume that $P_1=(-1,1)$, $P_2=O=(0,0)$
and $P_3=(1,1)$ and $\ell_1=\{y+2x+1=0\},\, \ell_3=\{y-2x+1=0\}$. 
Then the conic
must be defined by $y-x^2=0$. Thus $\ell_2=\{y=0\}$.
Here we used the next easy lemma.
\begin{Lemma}
Let $(C,\{P_1,P_2,P_3\})$ and $(C',\{P_1',P_2',P_3'\})$
two smooth conics with three points on the respective conic.
Then there are isomorphic by  an action  of a matrix $ A\in PGL(3,\bfC)$.
\end{Lemma}

Thus we need to have the equations
\[
 (\star):\begin{cases}
& f_{31}(P_i)=f_{32}(P_i)=0, \,  i=1,2,3\\
&(f_{3j,x}-\,2\,f_{3j,y})(P_1)=0,\,(f_{3j,x}+\,2\,f_{3j,y})(P_3)=0,\,
f_{3j,y}(O)=0, \, j=1,2\end{cases}
\]
The last condition says that two cubic are tangent to $y=x^2$ at $P_1,P_3$.
Let $R(x)$ and $S(y)$ be the resultant of $f_{31}$ and $f_{32}$ in
$x$-variable
and $y$-variable respectively. The above equality implies that 
$(x^2-1)^2x^2\,|\,R(x)$ and $y^2(y^2-1)^2\,|\,S(y)$.
Elliminating coefficients using these equalities, we consider the
further condition for $P_1,P_2,P_3$ to be $A_5$-singularities.
This is given by the condition 
$x^3(x^2-1)^3|R(x)$ and $y^3(y^2-1)^3|S(y)$.
At the end of calculation, we find that there are no such $f_{31},\, f_{32}$
which corresponds to a reducible sextics.

We consider the case (b).
Assume that $C_2$ is a union of two lines $L_1,L_2$. Then we can
see that $L_i$ are tangent to the cubic $C_3$ and the intersection
$L_1\cap L_2$ is also on $C_3$ so that this makes the third $A_5$.
In this case, we may assume that $P_1,P_2,P_3$ be as above but 
$\ell_1$ is $y=-x$ and $\ell_3$ is $y=x$.
The $\star$ should be replaced by
\[
 (\star):\begin{cases}
& f_{31}(P_i)=f_{32}(P_i)=0, \,  i=1,2,3\\
&(f_{3j,x}\,-\,f_{3j,y})(P_1)=0,\,(f_{3j,x}\,+\,f_{3j,y})(P_3)=0, \, j=1,2\end{cases}
\]
Then we consider the $A_5$-condition to see that there exists no such
sextics.

The case $\Si_{in}=[A_{11},A_5]$, we take $A_5$ at $P_1$ and $A_{11}$ at
$O$
and 
the $(\star)$-condition is repaleced by

\[
 (\star):\begin{cases}
& f_{31}(P_i)=f_{32}(P_i)=0, \,  i=1,2,
\quad (f_{3j,x}\,-\,2f_{3j,y})(P_1)=0,\, j=1,2\\
&x^4\,|\,f_{3j}(x,x^2), \, j=1,2\end{cases}
\]
The last condition says that the intersection multiplicity of each cubic
and the conic $y-x^2=0$  at $O$ is 4.
The reason that we have chosen the conic $y-x^2=0$ is to make the last
condition to be easier to be used.

The case $\Si_{in}=[A_{17}]$, we take $A_{17}$ at $P_2$, and the the
torus type condition is
\[
 (\star):
 f_{31}(0,0)=f_{32}(0,0)=0, \,
x^6\,|\,f_{3j}(x,x^2),\, j=1,2
\]
In any cases, one conclude that there does not exist any
solution
which corresponds to a reduced sextics.

\bibliographystyle{abbrv}
\def\cprime{$'$} \def\cprime{$'$}

\end{document}